\documentclass[reqno,12pt]{amsart}
\usepackage{amsfonts,amsmath,amssymb,amsthm,mathrsfs,array}
\usepackage[colorlinks, backref=section]{hyperref}
\usepackage[margin=2.4cm]{geometry}

\newtheorem{thm}{Theorem}[section]
\newtheorem{lemma}[thm]{Lemma}
\newtheorem{cor}[thm]{Corollary}
\newtheorem{prop}[thm]{Proposition}
\theoremstyle{definition}
\newtheorem{Def}[thm]{Definition}
\newtheorem{remark}[thm]{Remark}

\def\H{\mathcal{H}}

\def\A{\mathcal{A}}
\def\B{\mathcal{B}}
\def\R{\mathbb{R}}
\def\N{\mathbb{N}}
\def\Dt{\partial_t}
\def\Dtt{\partial_{tt}}
\def\Dx{\Delta}
\def\Nx{\nabla}

\def\ga{\gamma}
\def\al{\alpha}

\def\de{\delta}
\def\lm{\lambda}
\def\sg{\sigma}
\def\Om{\Omega}
\def\vp{\varphi}
\def\wto{\rightharpoonup}
\def\td{{\rm d}}
\def\({\left(}
\def\){\right)}
\def\le{\leqslant}
\def\ge{\geqslant}
\def\la{\langle}
\def\ra{\rangle}
\numberwithin{equation}{section}

\begin{document}

\title[Wave equations with low-regularity time-dependent forcing]{Well-posedness and pullback attractors for weakly damped wave equations
  with time-dependent low-regularity forcing}

\author{Cuncai Liu}
\email{liucc@jsut.edu.cn}
\address{School of Mathematics and Physics, Jiangsu University of Technology, Changzhou, 213001, China}

\author{Fengjuan Meng$^{*}$}
\email{fjmeng@jsut.edu.cn}
\address{School of Mathematics and Physics, Jiangsu University of Technology, Changzhou, 213001, China}

\author{Chang Zhang}
\email{chzhnju@126.com}
\address{School of Mathematics and Physics, Jiangsu University of Technology, Changzhou, 213001, China}
\author{Chunxiang Zhao}
\email{zhaocxmath@163.com}
\address{School of Mathematical Sciences, Jiangsu University,
Zhenjiang, 212013, China}
\thanks{*Corresponding author.\\
  This work was supported by the NSFC (11701230, 11801227, 11801228, 12026431), QingLan Project of Jiangsu Province and Zhongwu Young Innovative Talents Support Program.}

\begin{abstract}
  We consider the weakly damped semilinear wave equation
  $\Dtt u+\ga\Dt u-\Dx u+f(u)=g(t)$ on a smooth bounded domain
  $\Om\subset\R^n$, $n\ge3$, with Dirichlet boundary condition.
  The time-dependent force $g(t)$ belongs to the uniformly local
  space $H^1_b\(\R;H^{-1}(\Om)\)$; thus, its spatial regularity is
  only $H^{-1}(\Om)$, rather than the conventionally assumed $L^2(\Om)$.
  The nonlinearity is assumed to have growth order $p$, in the sense that
  $|f'(s)|\le C\(1+|s|^{p-1}\)$. For $1<p\le p_\Om(n)$, where
  $p_\Om(3)=25/7$ and $p_\Om(n)=1+2n/(n^2-3n+2)$ if $n\ge4$.
  We prove unconditional uniqueness of energy-class weak solutions without
  assuming any additional space-time regularity, by a duality argument.
  We then establish existence and Lipschitz continuous dependence on the
  initial data. The additional regularity is obtained for $u-\Theta$,
  where $\Theta$ solves the same equation with the same force and vanishing
  initial data. Finally, the resulting process possesses a pullback attractor
  in $\H=H^1_0(\Om)\times L^2(\Om)$, which is bounded in the past.
\end{abstract}


\maketitle
\setcounter{tocdepth}{1}

\section{Introduction}\label{sec-intro}

Let $\Om\subset\R^n$, $n\ge3$, be a bounded domain with smooth boundary. We consider the non-autonomous weakly damped wave equation
\begin{align}\label{eq1}
  \begin{cases}
    \Dtt u+\ga\Dt u-\Dx u+f(u)=g(t),             & (t,x)\in(\tau,+\infty)\times\Om,         \\
    u(\tau,x)=u_0(x),\quad \Dt u(\tau,x)=u_1(x), & x\in\Om,                                 \\
    u(t,x)=0,                                    & (t,x)\in(\tau,+\infty)\times\partial\Om,
  \end{cases}
\end{align}
where $\ga>0$ is a fixed dissipation rate. The phase space is the energy space
\[
  \H:=H^1_0(\Om)\times L^2(\Om),\qquad \xi_u(t):=\(u(t),\Dt u(t)\).
\]
We write $\|\cdot\|$ and $\la\cdot,\cdot\ra$ for the norm and scalar product of $L^2(\Om)$, $\|\cdot\|_q$ for the norm of $L^q(\Om)$, and $\lm_1>0$ for the first eigenvalue of $-\Dx$
with Dirichlet boundary condition. Throughout,
\begin{equation}\label{crit}
  2^*:=\frac{2n}{n-2},\qquad p_c:=2^*-1=\frac{n+2}{n-2}
\end{equation}
denote the Sobolev exponent of $H^1_0(\Om)$ and the associated critical growth exponent.

\subsection{Assumptions}
We impose the following conditions on the nonlinearity $f$ and on the external force $g$.
\begin{itemize}
  \item[({\bf F})] $f\in C^1(\R)$, $f(0)=0$, and $f$ satisfies the growth condition
    \begin{equation}\label{F1}\tag{F1}
      |f'(s)|\le C\(1+|s|^{p-1}\),\qquad 1< p\le p_\Om(n),
    \end{equation}
    where
    \[
      p_\Om(n):=
      \begin{cases}
        \dfrac{25}{7},          & n=3,   \\[2mm]
        1+\dfrac{2n}{n^2-3n+2}, & n\ge4,
      \end{cases}
    \]
    together with the dissipativity condition
    \begin{equation}\label{F2}\tag{F2}
      \begin{aligned}
         & f(s)s\ge-\lm s^2-C_f,\qquad F(s):=\int_0^sf(\varsigma)\,\td\varsigma
        \ge-\frac{\lm}{2}s^2-C_f,                                               \\
         & f(s)s-F(s)\ge-\frac{\lm}{2} s^2-C_f,
      \end{aligned}
    \end{equation}
    for some constants $C_f\ge0$ and $0\le\lm<\lm_1$, where $\lm_1$ is the first eigenvalue of $-\Delta$ with Dirichlet boundary condition on $\Omega$.

  \item[({\bf G})] The external force satisfies
    \begin{equation}\label{G1}\tag{G1}
      g\in H^1_b\(\R;H^{-1}(\Om)\),
    \end{equation}
    where $H^1_b$ is the uniformly local Sobolev space of Section \ref{ssec-forces}.
\end{itemize}

The threshold $p_\Om(n)$ is dictated by the Strichartz estimates available on a bounded domain; its derivation is given in Appendix \ref{app-opt}. Since $p_\Om(n)<p_c$ in every
dimension, integrating \eqref{F1} gives
\begin{equation}\label{fbounds}
  |f(s)|\le C\(1+|s|^{p}\),\qquad |F(s)|\le C\(1+|s|^{p+1}\),\qquad p+1<2^*,
\end{equation}
so that, by \eqref{F2}, the energy functional of \eqref{eq1} is well defined and bounded from below on $\H$.

Hypothesis \eqref{G1} says that the force has only the spatial regularity of $H^{-1}(\Om)$, one derivative below the classical assumption $g\in L^2(\Om)$. It has two consequences, which
pull in opposite directions.
\begin{itemize}
  \item It does \emph{not} give $g\in L^1_{loc}(\R;L^2(\Om))$. This single fact governs much of the analysis: it rules out any direct application of a Strichartz estimate to a solution of
        \eqref{eq1}, and it is the main difficulty created by a low-regularity forcing.
  \item The time derivative in \eqref{G1} is what makes an energy estimate possible at all. The pairing $\la g(t),\Dt u(t)\ra$ is meaningless for $g(t)\in H^{-1}(\Om)$ and
        $\Dt u(t)\in L^2(\Om)$, and it has to be replaced by $\frac{\td}{\td t}\la g,u\ra-\la\Dt g,u\ra$.
\end{itemize}

The purpose of this paper is to prove existence and uniqueness of weak solutions of \eqref{eq1} in the energy space $\H$ under the low-regularity assumption \eqref{G1}, and to establish the existence of the pullback attractor of the process generated by these solutions.

\subsection{Time-independent low-regularity forcing}

When the force does not depend on time, the difficulty created by $g\in H^{-1}$ can be removed by subtracting a stationary solution. In our previous works \cite{lcc-tmna,lcc-jde,lcc-aml,lcc26}, we used
this to obtain well-posedness in the energy space for $n=3$. We describe the mechanism for the equation on $\R^3$, where a mass term $m u$, $m>0$, is added to \eqref{eq1} so that $-\Dx+m$ is
coercive on $H^1(\R^3)$.

\emph{The cubic case $p=3$.} Here $|f(u)|\lesssim|u|+|u|^3$, and $H^1(\R^3)\hookrightarrow L^6(\R^3)$ gives $f(u)\in L^2(\R^3)$. Pairing the equation with $\Dt u$
therefore yields a genuine energy identity, and uniqueness follows from standard $L^2$ estimates for the difference of two solutions \cite{lcc-aml}.

\emph{The range $3<p<4$.} For $p>3$ the term $f(u)$ no longer lies in $L^2(\R^3)$, and the argument above fails. Let $h\in H^1(\R^3)$ be a solution of the stationary problem
$-\Dx h+mh+f(h)=g$. Although $f(u)$ itself is not square integrable, the difference of the nonlinearities is controlled by H\"older's inequality:
\begin{equation}\label{u-v}
  \|f(u)-f(h)\|_{L^2(\R^3)}
  \le C\(\|u\|_{L^6}^{p-1}+\|h\|_{L^6}^{p-1}\)\|u-h\|_{L^{\frac{6}{4-p}}(\R^3)},
  \qquad 3<p<4,
\end{equation}
for the superlinear part of $f$; the linear part is harmless. Hence, if $u-h$ has the higher integrability $u-h\in L^{6/(4-p)}(\R^3)$, then $f(u)-f(h)\in L^2(\R^3)$, and the equation for
$u-h$ admits energy estimates. This higher integrability comes from Strichartz estimates. In \cite{lcc-tmna} we introduced \emph{translationally regular solutions}, namely those $u$ for
which
\[
  u-h\in L^{\frac{2}{p-3}}_t L^{\frac{6}{4-p}}_x(\R^3).
\]
For such solutions the nonlinear difference acts as an $L^2$ source, and an iteration argument gives global well-posedness and continuous dependence on the data. The exponent pair above is
admissible on $\R^3$ exactly when $p<4$. On a bounded domain the smaller admissible range of \cite{blair09} restricts the same argument to $p\le\frac{25}{7}$ \cite{lcc-jde}.

\emph{The range $4\le p<5$.} For $p\ge4$ the exponent $\frac{6}{4-p}$ is no longer a positive finite exponent, and \eqref{u-v} breaks down, so a different idea is needed.
In \cite{lcc26}, assuming in addition $f\in C^2$, we expand around the stationary solution. The perturbation $\de(t):=u(t)-h$ satisfies
\begin{equation}\label{delta_intro}
  \Dtt\de+\ga\Dt\de-\Dx\de+m\de+f(u)-f(h)=0,
\end{equation}
and
\begin{equation}\label{fu-fv}
  f(u)-f(h)=f'(h)\,\de+\de^2\int_0^1(1-\tau)f''(h+\tau\de)\,\td\tau .
\end{equation}
Setting
\[
  V(x):=f'\(h(x)\),\qquad
  N(t,x):=-\int_0^1(1-\tau)f''\(h(x)+\tau\de(t,x)\)\,\td\tau,
\]
equation \eqref{delta_intro} becomes a wave equation with a time-independent potential,
\begin{equation}\label{delta_V}
  \Dtt\de+\ga\Dt\de-\Dx\de+V(x)\de=N(t,x)\,\de^2-m\de .
\end{equation}
Strichartz estimates for \eqref{delta_V}, obtained in the Besov spaces associated with the Schr\"odinger operator $-\Dx+V$, then give global well-posedness of weak solutions through
an iteration on the Besov scale.

In the above cases the decisive step is the same: the low-regularity forcing is absorbed into a stationary solution $h$. This requires $g$ to be independent of time, and it fails for a
time-dependent force $g(t)$.

\subsection{The time-dependent forcing}
In the present time-dependent forcing case, equation \eqref{eq1} generally has no stationary solution. Thus, unlike the time-independent setting treated in \cite{lcc-tmna,lcc-jde,lcc26}, the rough force cannot be absorbed into a stationary solution. To absorb the force in the time-dependent case, a suitable replacement for $h$ is the solution $\Theta$ of \eqref{eq1} with vanishing initial data, namely the solution of
\begin{align}\label{eqTheta}
  \partial_{t}^2 \Theta+\ga \Dt \Theta-\Dx \Theta+f(\Theta)=g(t), \quad \Theta(\tau,x)=0,\quad \Dt\Theta(\tau,x)=0,  x\in \Omega.
\end{align}
For any weak solution $u$ of \eqref{eq1}, the difference $u-\Theta$ solves an equation in which $g$ cancels identically, and whose right-hand side again has the difference structure $f(\Theta)-f(u)$. Two new difficulties arise.
\begin{itemize}
  \item In \cite{lcc-tmna, lcc-jde}, we introduced the \emph{translationally regular solutions},
        namely the weak solutions $u$ satisfying extra integrability $u-h\in L^{\frac{2}{p-3}}_t L^{\frac{6}{4-p}}_x(\R^3)$. For time-dependent forcing, one might similarly define a subclass of weak
        solutions by requiring $u-\Theta\in L^q_t L^r_x(\Om)$ for some $q\ge 2$ and $r>2^*$.
        For this definition to be independent of the choice of $\Theta$, one would need to show that
        $\Theta_1-\Theta_2\in L^q_t L^r_x(\Om)$ for any two energy-class weak solutions $\Theta_1$ and $\Theta_2$ of \eqref{eqTheta}, without assuming Strichartz space-time regularity. This amounts to
        proving unconditional uniqueness in the energy class. Once uniqueness is established, there is no need to single out a subclass since the weak solution itself is unique. We therefore first prove unconditional uniqueness by adapting the duality method of \cite{glassey82} to a bounded domain. As a byproduct, the \emph{translationally regular solutions} of \cite{lcc-tmna, lcc-jde} coincide with weak solutions.
  \item The potential $f'(\Theta)$ now depends on time and has only the integrability supplied by the energy bound. The Besov-space approach used for $4\le p<5$ is then unavailable,
        and only the H\"older mechanism of \eqref{u-v} remains. On a bounded domain, optimising this mechanism over all admissible Strichartz exponents gives exactly the threshold
        $p_\Om(n)$; see Appendix \ref{app-opt}.
\end{itemize}

\subsection{Uniqueness of weak solutions}

The natural benchmark for the admissible range of $p$ is the undamped, unforced wave equation
\begin{equation}\label{nlw}
  \Dtt u-\Dx u+f(u)=0\qquad\text{on }\Om\text{ or on }\R^n .
\end{equation}
For energy data, treating $f(u)$ as a source term, the maximal source norm available in the energy class is $L^1(I;L^2_x)$. Requiring $f(u)\in L^2_x$ means $u\in L^{2p}_x$, whereas the
energy bound gives only $u\in L^{2^*}_x$. Hence uniqueness of weak solutions follows from a direct energy argument as long as
\begin{equation}\label{naive}
  p\le\frac{2^*}{2}=\frac{n}{n-2},
\end{equation}
and $\frac{n}{n-2}$ is the first obstruction.

On $\R^n$, Glassey and Tsutsumi \cite{glassey82} proved unconditional uniqueness in the locally finite energy class for $1<p\le p_1(n):=\frac{n^2+3n-2}{n^2-n-2}$ by a duality method. Note that $\frac{n}{n-2}<p_1(n)<p_c$. Through an iteration of regularity, Ginibre and Velo \cite{Ginibre} obtained unconditional uniqueness in the whole subcritical range $p<p_c$, $n\ge3$.The critical exponent $p=p_c$ is considerably more delicate, and has been reached only in higher dimensions. Planchon \cite{planchon03} established unconditional uniqueness in $C_t(\dot H^1)$ for
$n=4,5,6$. For $n\ge6$ the same conclusion was obtained by Bulut et al. \cite{bulut}, in these dimensions the nonlinearity fails to be $C^2$, so the second-order Taylor expansion available in the low-dimensional case is not available. In three
dimensions, i.e. at the quintic exponent $p=5$, unconditional uniqueness in the energy space remains open. Only conditional results are known there: uniqueness holds in the Strichartz
class \cite{ss}, and among solutions satisfying a strong energy inequality \cite{struwe99} or local energy inequality \cite{masmoudi}.

On a bounded domain the situation is more delicate: Burq, Lebeau and Planchon \cite{burq} proved global well-posedness and uniqueness for the defocusing energy-critical equation in the critical Strichartz class $L^5_tL^{10}_x$.

\subsection{Long-term behavior}
The long-term behavior of dissipative wave equations is classical, and we refer to the monographs
\cite{babin92,haraux87,temam97} for the general theory of global attractors. For non-autonomous
The theory of \emph{uniform} attractors for a family of processes indexed by the hull of the force
is due to Chepyzhov and Vishik \cite{chueshov02}. Producing a strong uniform attractor requires the
hull to be compact in a strong topology, which in turn requires the force to be translation
compact. Zelik \cite{zelik15} weakened this requirement by introducing the classes of \emph{space regular} and \emph{time regular} forces, either of which separately suffices to build a strong
uniform attractor for an $L^2(\Om)$-valued force, while their conjunction is equivalent to
translation compactness Theorem 2.8 of \cite{zelik15}.

The alternative notion, which we adopt here, is the \emph{pullback} attractor of the single process
generated by the fixed force $g$. Its existence does not require $g$ to belong to a compact family of translates, and hence no space or time regularity of $g$ beyond \eqref{G1} is needed: pullback asymptotic compactness is obtained by the energy method of Ball \cite{ball04} and of Moise, Rosa and Wang \cite{mrw98}.

\subsection{Results}
The contribution of the paper can be summarised in three points. First, we prove unconditional uniqueness in the
energy class for the time-dependent $H^{-1}$-valued forcing problem. This is stronger than a uniqueness statement
restricted to a prescribed Strichartz class, and it is the step that removes the apparent circularity caused by the
rough forced component $\Theta$. Second, once uniqueness is available, we recover the largest exponent admitted by the
present bounded-domain Strichartz--duality scheme, namely $p\le p_\Om(n)$, and we obtain extra space-time regularity
for $u-\Theta$ rather than for $u$ itself, which is the most the hypothesis \eqref{G1} can support. Third, the
resulting well-posed process admits a pullback attractor under the single assumption that $g$ is translation bounded
and $H^1$ in time. The argument therefore combines the bounded-domain Strichartz estimates of \cite{blair09}, a
regularised duality argument, and the energy method of \cite{ball04,mrw98}, while avoiding any
translation-compactness assumption of the kind used in the uniform-attractor theory of \cite{chueshov02,zelik15}.

The order in which the results of Section \ref{sec-wp} are proved is not interchangeable, a point worth stating
explicitly because the natural arrangement is circular. Uniqueness comes first and is \emph{unconditional}, in the
energy class and with no space-time integrability assumed; only then is the auxiliary solution $\Theta$ of
\eqref{eqTheta} well defined; and only then can existence be proved together with the extra regularity of
$u-\Theta$.

\begin{thm}[Unconditional uniqueness; Theorem \ref{thm-unique}]\label{unique}
  Assume \eqref{F1}. Then, for every $T>\tau$, problem \eqref{eq1} admits at most one weak solution in the energy class of Definition \ref{def-weak}.
\end{thm}
We next obtain existence together with the extra regularity.
\begin{thm}[Existence, extra regularity, Lipschitz dependence; Theorems \ref{thm-exist} and
    \ref{thm-lip}]\label{thmC}
  Assume \eqref{F1}, \eqref{F2} and \eqref{G1}, let
  \begin{equation}\label{optpair}
    (q,r)=
    \begin{cases}
      \Big(\dfrac{7}{2},\ 14\Big),         & n=3,   \\[3mm]
      \Big(n-1,\ \dfrac{2(n-1)}{n-3}\Big), & n\ge4,
    \end{cases}
  \end{equation} and let $\Theta$ be the unique weak solution of \eqref{eq1} with vanishing data at $t=\tau$. Then
  for every $\xi_\tau\in\H$ problem \eqref{eq1} has a unique weak solution, and
  \[
    \|\xi_u(t)\|_{\H}+\|u-\Theta\|_{L^{q}(t,t+1;L^{r}(\Om))}
    \le Q\(\|\xi_\tau\|_{\H}+\|g\|_{H^1_b}\)e^{-\al(t-\tau)}+Q\(\|g\|_{H^1_b}\),\qquad t\ge\tau,
  \]
  for a monotone $Q$ independent of $\tau,t,\xi_\tau$. Moreover any two solutions satisfy
  \[
    \|\xi_{u^1}(t)-\xi_{u^2}(t)\|_{\H}+\|u^1-u^2\|_{L^{q}(\tau,t;L^{r}(\Om))}
    \le Ce^{K(t-\tau)}\|\xi_{u^1}(\tau)-\xi_{u^2}(\tau)\|_{\H}.
  \]
\end{thm}

The second term on the left involves $u-\Theta$ and not $u$: under \eqref{G1} we cannot assert that an individual solution lies in $L^{q}_{loc}(L^{r})$. The force must first be removed by subtraction, which is possible only because $\Theta$ solves the equation with the \emph{same} force.

\begin{thm}[Pullback attractor; Theorem \ref{thm-pullback}]\label{thmD}
  Assume \eqref{F1}, \eqref{F2} and \eqref{G1}. Then the process $U(t,\tau):\H\to\H$,
  $U(t,\tau)\xi_\tau:=\xi_u(t)$, generated by the weak solutions of \eqref{eq1}, possesses a pullback attractor $\A(\cdot)$, bounded in the past, with the representation
  \[
    \A(t)=\big\{u(t):u\ \text{is a complete weak solution of \eqref{eq1} that is bounded backwards from }t\big\}.
  \]
\end{thm}

Besides well-posedness and the dissipative estimate, the only ingredient of Theorem \ref{thmD}
is the pullback asymptotic compactness of $U(\cdot,\cdot)$. Fix $t\in\R$ and $s_k\to-\infty$, and
let $u_k$ be solutions with data bounded at time $s_k$. The dissipative estimate bounds $\xi_{u_k}$
in $\H$ uniformly on $[s_k,t]$, and the Aubin--Lions--Simon theorem then yields a limit trajectory
$u$ on $(-\infty,t]$ with $\xi_{u_k}(t)\wto\xi_u(t)$ weakly in $\H$; it yields strong convergence
only in $H^{1-\de}(\Om)\times H^{-\de}(\Om)$. Strong convergence in $\H$ is recovered by the
energy method of Ball \cite{ball04} and of Moise, Rosa and Wang \cite{mrw98}: we apply the
exponentially weighted energy identity of Lemma \ref{lem-expenergy} to each $u_k$ on $[s_k,t]$ and
to $u$ on $(-\infty,t]$. The weight $e^{\ga s}$ eliminates the dissipation integral, which is only
lower semicontinuous under weak convergence, and causes the contribution of the initial time $s_k$
to vanish as $k\to\infty$. This yields $\|\xi_{u_k}(t)\|_{\H}\to\|\xi_u(t)\|_{\H}$. Since only the single force $g$ enters this argument, the hull $\mathcal H(g)$ of Chepyzhov and Vishik \cite{chueshov02}, and the classes of space regular and time regular forces of Zelik \cite{zelik15}
which replace translation compactness in the \emph{uniform} attractor theory of a family of
processes, play no role here.

Throughout, $C$ denotes a generic positive constant that may change from line to line, whereas $C_0,C_1,\dots$ are fixed once they have been introduced.

\section{Preliminaries}\label{sec-prelim}

\subsection{Fractional spaces}
Let $\{e_i\}_{i\in\N}$ be the orthonormal basis of $L^2(\Om)$ consisting of Dirichlet
eigenfunctions of $-\Dx$, with $0<\lm_1\le\lm_2\le\cdots\to+\infty$. For $s\in\R$ set
\[
  H^{s}:=D\((-\Dx)^{s/2}\)
  =\Big\{w=\sum_{i\ge1}a_ie_i:\ \|w\|_{H^s}^2
  :=\sum_{i\ge1}\lm_i^{s}a_i^2<\infty\Big\},
\]
so that $H^0=L^2(\Om)$, $H^1=H^1_0(\Om)$ and $H^{-1}=\(H^1_0(\Om)\)'$. Let $P_N$ be the orthogonal
projection onto $\mathrm{span}\{e_1,\dots,e_N\}$; since the $e_i$ are eigenfunctions, $P_N$ is orthogonal in every $H^{s}$ at once. We shall use freely the interpolation inequality
$\|w\|_{H^{\theta s_0+(1-\theta)s_1}}
  \le\|w\|^{\theta}_{H^{s_0}}\|w\|^{1-\theta}_{H^{s_1}}$ and the embedding
$H^1_0(\Om)\hookrightarrow L^{2^*}(\Om)$.

\subsection{Uniformly local spaces}\label{ssec-forces}
For a reflexive Banach space $V$ and $1<\theta<\infty$ put
\[
  L^{\theta}_b(\R;V):=\Big\{h\in L^{\theta}_{loc}(\R;V):\
  \|h\|_{L^\theta_b(\R;V)}:=\sup_{t\in\R}\|h\|_{L^\theta(t,t+1;V)}<\infty\Big\},
\]
and define $H^k_b(\R;V)$ analogously, with $\sum_{j\le k}\|\Dt^jh(s)\|_V^2$ under the
integral. Elements of $L^\theta_b$ are called \emph{translation bounded}, and hypothesis \eqref{G1} says precisely that $g\in H^1_b(\R;H^{-1}(\Om))$. The continuous embedding
$H^1(t,t+1;H^{-1})\hookrightarrow C([t,t+1];H^{-1})$, whose constant does not depend on $t$, yields
\begin{equation}\label{gsup}
  \sup_{t\in\R}\|g(t)\|_{H^{-1}}\le C_*\|g\|_{H^1_b(\R;H^{-1})}<\infty .
\end{equation}

\begin{remark}\label{rem-no-hull}
  In the uniform-attractor theory of Chepyzhov and Vishik \cite{chueshov02} one works with the
  entire hull $\mathcal H(g)$ of translates of $g$, that is, the closure of $\{T_\ell g:\ell\in\R\}$,
  $T_\ell g(t)=g(t+\ell)$, taken in the local weak topology $L^2_{\rm loc}(\R\times\Om;\mathcal D'(\Om))$,
  and the asymptotic compactness required for a \emph{strong} uniform attractor must hold \emph{uniformly} over this hull.
  Since a generic $h\in\mathcal H(g)$ need not be any more regular than $g$ itself, Zelik \cite{zelik15}
  introduced the classes of space regular and time regular forces precisely in order to control this uniformity when $g$ is $L^2(\Om)$-valued.

  The pullback attractor of Section \ref{sec-attr} is built instead from the single process generated by the fixed
  force $g$: the only limits taken are along sequences $s_k\to-\infty$ with the \emph{same} force $g$, never along
  a sequence of distinct translates or approximants of $g$. Consequently none of the apparatus above --- the hull,
  its weak compactness, or any class of forces substituting for translation compactness --- is needed, and
  \eqref{G1} alone is the only hypothesis used in Section \ref{sec-attr}.
\end{remark}


\subsection{Strichartz estimates on a bounded domain}
Consider the linear problem
\begin{equation}\label{linwave}
  \Dtt w-\Dx w=\mathcal F(t)\ \text{ in }\Om\times\R,\qquad w|_{\partial\Om}=0,\qquad
  (w, \Dt w)|_{t=0}=(w_0,w_1).
\end{equation}

\begin{Def}\label{def-adm}
  Let $n\ge3$. A triple $(q,r,\sg)$ with $2<q\le\infty$, $2\le r<\infty$ is
  \emph{$\Om$-admissible} if
  \begin{equation}\label{adm-om}
    \frac1q+\frac nr=\frac n2-\sg
    \qquad\text{and}\qquad
    \begin{cases}
      \dfrac3q+\dfrac{n-1}{r}\le\dfrac{n-1}{2}, & n=3,   \\[4mm]
      \dfrac1q+\dfrac1r\le\dfrac12,             & n\ge4,
    \end{cases}
  \end{equation}
  and \emph{$\R^n$-admissible} if the scaling identity holds together with
  $\frac2q+\frac{n-1}{r}\le\frac{n-1}{2}$, which is the standard condition on $\R^n$.
\end{Def}
Blair, Smith and Sogge \cite{blair09} established the following Strichartz estimates for the wave equation on a manifold with boundary.
\begin{thm}\label{thm-BSS}
  Let $\Om$ be a compact Riemannian manifold with boundary of dimension $n$, let $(q,r,\sg)$ be
  $\Om$-admissible, and impose either the Dirichlet or the Neumann boundary condition. Then every solution of \eqref{linwave} with $\mathcal F=0$ satisfies, for every $T>0$,
  \[
    \|w\|_{L^q(-T,T;L^r(\Om))}\le C_T\(\|w_0\|_{H^{\sg}}+\|w_1\|_{H^{\sg-1}}\),
  \]
  and if $(\tilde q,\tilde r,1-\sg)$ is $\Om$-admissible as well, then
  \begin{equation}\label{BSS-inhom}
    \|w\|_{L^q(-T,T;L^r(\Om))}
    \le C_T\(\|w_0\|_{H^{\sg}}+\|w_1\|_{H^{\sg-1}}
    +\|\mathcal F\|_{L^{\tilde q'}(-T,T;L^{\tilde r'}(\Om))}\)
  \end{equation}
  holds, where $\tilde q'$ and $\tilde r'$ are the H\"older conjugates of $\tilde q$ and $\tilde r$. The constant $C_T$ depends on $T,\Om,n,q,r,\sg$.
\end{thm}
We also need the corresponding statement for the damped equation
\begin{equation}\label{dampedlin}
  \Dtt w+\ga\Dt w-\Dx w=\mathcal F\ \text{ in }\Om\times\R,\qquad w|_{\partial\Om}=0,
  \qquad (w,\Dt w)|_{t=0}=(w_0,w_1).
\end{equation}
\begin{lemma}\label{lem-damped}
  Let $(q,r)$ be as in \eqref{optpair} and let $I=[\tau,\tau+L]$ with $L\le1$. Then every
  solution of \eqref{dampedlin} with $(w_0,w_1)\in\H$ satisfies
  \begin{equation}\label{main-str}
    \|w\|_{L^{q}(I;L^{r}(\Om))}\le C_0\(\|w_0\|_{H^1_0}+\|w_1\|
    +\|\mathcal F\|_{L^1(I;L^2(\Om))}\),
  \end{equation}
  with $C_0=C_0(\ga,\Om,n)$ independent of $\tau$, $w_0$, $w_1$ and $\mathcal F$.
\end{lemma}

\begin{proof}
  Set $v(t):=e^{\ga(t-\tau)/2}w(t)$, so that $v|_{\partial\Om}=0$ and
  \begin{equation}\label{vshift}
    \Dtt v-\Dx v=\frac{\ga^2}{4}v+e^{\ga(t-\tau)/2}\mathcal F,\qquad
    v(\tau)=w_0,\quad \Dt v(\tau)=\tfrac\ga2 w_0+w_1 .
  \end{equation}
  Since $L\le1$ the factor $e^{\ga(t-\tau)/2}$ is bounded by $e^{\ga/2}$ on $I$ and is
  harmless in every norm below. Applying \eqref{BSS-inhom} with $\sg=1$,
  $(\tilde q,\tilde r)=(\infty,2)$ and the whole right-hand side of \eqref{vshift} as source,
  \[
    \|v\|_{L^q(I;L^r)}\le C\(\|w_0\|_{H^1_0}+\|w_1\|
    +\big\|\tfrac{\ga^2}{4}v\big\|_{L^{1}(I;L^{2})}
    +\|\mathcal F\|_{L^{1}(I;L^{2})}\).
  \]
  The term containing $v$ is \emph{not} absorbed by smallness, which would be circular, but
  estimated independently of the Strichartz norm: applying the energy estimate to \eqref{dampedlin},
  we have
  \[
    \sup_{t\in I}\big(\|\Dt w(t)\|+\|\Nx w(t)\|\big)
    \le C\(\|w_0\|_{H^1_0}+\|w_1\|+\|\mathcal F\|_{L^1(I;L^2)}\),
  \]
  whence
  \[
    \|v\|_{L^1(I;L^2)}\le Le^{\ga L/2}\sup_{I}\|w\|\le C\(\|w_0\|_{H^1_0}+\|w_1\|+\|\mathcal F\|_{L^1(I;L^2)}\).
  \]
  Returning to $w=e^{-\ga(t-\tau)/2}v$ gives \eqref{main-str}.
\end{proof}
\subsection{Pullback attractors}

For the convenience of the reader, we recall some standard definitions from theory of pullback
attractors, for more details see \cite{cara04,clr13} and references therein.

Let $X$ be a metric space and $\{U(t, \tau):X\to X \mid -\infty<\tau\le t<+\infty\}$ be
a process on $X$ satisfying
\begin{itemize}
  \item[(i)]  $U(\tau, \tau)=\mathrm{Id}$;
  \item[(ii)]  $U(t, \tau) = U(t, s)U(s,\tau )$, for all $\tau\le s\le t$.
\end{itemize}

\begin{Def}\label{def-pb-attraction}
  Let $\{U(t,\tau)\}_{t\ge\tau}$ be a process on a metric space $X$. Given $\A$ and $\B$
  subsets of $X$, we say that $\A$ \emph{pullback attracts} $\B$ at time $t\in\R$ if
  $$\lim_{s\to-\infty}\mathrm{dist}_{X}\big(U(t,s)\B, \A\big)=0,$$
  where $\mathrm{dist}_X$ denotes the Hausdorff semi-distance. $\A$ pullback attracts
  bounded sets at time $t$ if this holds for every bounded $\B\subset X$.
\end{Def}
\begin{Def}\label{def-invariant}
  A family of subsets $\{\A(t):t\in\R\}$ of $X$ is \emph{invariant} relative to the process
  $\{U(t,\tau)\}_{t\ge\tau}$ if $U(t,\tau)\A(\tau)=\A(t)$ for every $t\ge\tau$.
\end{Def}
\begin{Def}\label{def-pb-attractor}
  A family $\{\A(t):t\in\R\}$ of subsets of $X$ is the \emph{pullback attractor} of $U(\cdot,\cdot)$ if (i)
  $\A(\cdot)$ is invariant; (ii) $\A(t)$ is compact in $X$ for every $t$; (iii) $\A(\cdot)$ pullback attracts
  bounded subsets of $X$; and (iv) $\A(\cdot)$ is the minimal family of closed sets having property (iii).
\end{Def}

Two further notions are needed for the existence theorem.
\begin{Def}\label{def-pb-dissip}
  $U(\cdot,\cdot)$ is \emph{strongly pullback bounded dissipative} if there is a bounded set $\B_0\subset X$,
  independent of $t$, such that $\B_0$ pullback attracts bounded sets at time $\varsigma$ for every
  $\varsigma\in\R$.
\end{Def}
\begin{Def}\label{def-pb-ac}
  $U(\cdot,\cdot)$ is \emph{pullback asymptotically compact} if, for every $t\in\R$, every sequence
  $\tau_k\to-\infty$ and every bounded sequence $\{x_k\}\subset X$, the sequence $\{U(t,\tau_k)x_k\}$ has a
  convergent subsequence in $X$.
\end{Def}

\begin{thm}[Existence; Theorem 3.1 of \cite{clr13}]\label{thm-CLR-existence}
  If $U(\cdot,\cdot)$ is strongly pullback bounded dissipative and pullback asymptotically compact, then it
  possesses a pullback attractor $\A(\cdot)$, bounded in the past, given by $\A(t)=\omega(\B_0,t)$, the pullback
  $\omega$-limit set of $\B_0$ at time $t$.
\end{thm}

\begin{thm}[Structure of the attractor; Theorem 6.4 of \cite{clr13}]\label{thm-CLR-kernel}
  If the pullback attractor $\A(\cdot)$ is bounded in the past, then
  \[
    \A(t)=\big\{\varphi(t):\ \varphi\ \text{is a global solution with } \sup_{s\le t}\|\varphi(s)\|_{\H}<\infty\big\}.
  \]
\end{thm}
\section{Well-posedness}\label{sec-wp}
\subsection{Weak solution}\label{ssec-ws}
\begin{Def}\label{def-weak}
  Let $\xi_\tau=(u_0,u_1)\in\H$ and $T>\tau$. A function $u$ with
  \[
    u\in L^{\infty}\(\tau,T;H^1_0(\Om)\),\qquad \Dt u\in L^{\infty}\(\tau,T;L^2(\Om)\),\qquad
    u(\tau)=u_0,\ \Dt u(\tau)=u_1,
  \]
  is a \emph{weak solution} of \eqref{eq1} on $[\tau,T]$ if
  \[
    \int_\tau^T\!\!\int_\Om\big[-\Dt u\,\Dt\psi+\ga\,\Dt u\,\psi
      +\Nx u\cdot\Nx\psi+f(u)\psi\big]\,\td x\,\td t
    =\int_\tau^T\big\la g(t),\psi(t)\big\ra\,\td t
  \]
  for every $\psi\in C^\infty_0\((\tau,T)\times\Om\)$.
\end{Def}

Let $(\tilde q,\tilde r,1-\sg)$ be any admissible dual triple. Hypothesis \eqref{G1} says that
$g(t)\in H^{-1}(\Om)$, and since $\tilde r'\ge1$ while $H^{-1}(\Om)\not\subset L^1(\Om)$, in general
\[
  g\notin L^{\tilde q'}(I;L^{\tilde r'}(\Om))
\]
for every interval $I$. Consequently, under \eqref{G1} no Strichartz estimate can be applied to a solution of
\eqref{eq1} taken alone, in any dimension, for any admissible configuration, and for any value of $p$, however
small. The force must therefore be removed before any Strichartz estimate is invoked, and the only exact way to
remove it is to subtract a solution carrying the \emph{same} force. That solution must be well defined, which means
that uniqueness has to be established first, without recourse to any space-time regularity. Accordingly,
Section \ref{ssec-unique} proves unconditional uniqueness, Section \ref{ssec-exist} proves existence together with
the extra regularity, and Section \ref{ssec-lip} proves Lipschitz dependence.

\subsection{Unconditional uniqueness}\label{ssec-unique}

\begin{thm}\label{thm-unique}
  Assume \eqref{F1}. Let $\xi_\tau\in\H$ and $T>\tau$. Then problem \eqref{eq1} admits at most one weak solution on $[\tau,T]$.
\end{thm}

The proof is a duality argument applied to the equation for the difference
$w:=u-v$. The same method was used in \cite{glassey82} to achieve uniqueness in the case
$p\le p_1(n)$ on the whole space $\R^n$, and in \cite{masmoudi} to prove uniqueness for weak solutions satisfying a local energy inequality in the case $p=5$ and $n=3$. Let $u,v$ be two weak solutions with the same data and set $w:=u-v$. Subtracting the two equations, the force cancels identically:
\begin{equation}\label{w-eq}
  w_{tt}+\ga \Dt w-\Dx w+Vw=0,\qquad w(\tau)=\Dt w(\tau)=0,
\end{equation}
where
\[
  V:=\int_0^1f'\(su+(1-s)v\)\,\td s,\qquad Vw=f(u)-f(v).
\]
Set
\begin{equation}\label{mu}
  \mu:=\frac{2^*}{p-1},\qquad
  R:=\sup_{t\in[\tau,T]}\max\{\|\xi_u(t)\|_{\H},\,\|\xi_v(t)\|_{\H},\,1\}<\infty.
\end{equation}
By \eqref{F1} and $H^1_0(\Om)\hookrightarrow L^{2^*}(\Om)$, the potential $V$ satisfies
\[
  |V|\le C\(1+|u|^{p-1}+|v|^{p-1}\),
\]
and therefore, on $I=[\tau,T]$,
\begin{equation}\label{V-R}
  \|V\|_{L^\infty(I;L^{\mu}(\Om))}
  \le C_1\(1+\sup_{I}\|\xi_u\|_{\H}^{p-1}+\sup_{I}\|\xi_v\|_{\H}^{p-1}\)
  \le C_1\(1+2R^{p-1}\).
\end{equation}

We remove the damping from \eqref{w-eq} by setting
\begin{equation}\label{delta-def}
  \de(t):=e^{\ga(t-\tau)/2}w(t),\qquad m:=\frac{\ga^2}{4}.
\end{equation}
Since $w_{tt}+\ga\Dt w=e^{-\ga(t-\tau)/2}(\de_{tt}-m\de)$, equation \eqref{w-eq} becomes
\begin{equation}\label{delta-eq}
  \de_{tt}-\Dx\de+(V-m)\de=0,\qquad \de(\tau)=\Dt\de(\tau)=0,
\end{equation}
which contains no first-order time derivative. The function $\de$ inherits the regularity of $w$, so
$\de\in L^{\infty}(\tau,T;H^1_0(\Om))$ and $\Dt\de\in L^{\infty}(\tau,T;L^2(\Om))$. Moreover,
\[
  \frac1\mu+\frac1{2^*}=\frac{p-1}{2^*}+\frac1{2^*}=\frac{p}{2^*},
\]
so H\"older's inequality yields
\[
  \|V\de\|_{L^{2^*/p}(\Om)}\le\|V\|_{L^{\mu}(\Om)}\|\de\|_{L^{2^*}(\Om)}.
\]
Combining with $L^{2^*/p}(\Om) \subset  H^{-1}(\Om)$, since $\frac{2^*}{p}\ge(2^*)'$, yields
\begin{equation}\label{Vdelta}
  V\de\in L^\infty\(\tau,T;L^{2^*/p}(\Om)\)\subset L^\infty\(\tau,T;H^{-1}(\Om)\),
\end{equation}
Hence \eqref{delta-eq} holds in $L^{\infty}(\tau,T;H^{-1}(\Om))$ and
$\de_{tt}\in L^{\infty}(\tau,T;H^{-1}(\Om))$.

\subsubsection*{The regularised dual problem}
Fix $F\in C^\infty_0\((\tau,T)\times\Om\)$. We test \eqref{delta-eq} against the solution of the backward problem
with source $F$,
\begin{equation}\label{dual}
  \vp_{tt}-\Dx\vp+(V-m)\vp=F\ \text{ in }(\tau,T)\times\Om,\quad
  \vp|_{\partial\Om}=0,\quad \vp(T)=\Dt\vp(T)=0 .
\end{equation}
Solving \eqref{dual} directly would require an approximation argument followed by passage to the limit, so we
instead work with mollified equations from the outset. Since $f'$ is merely continuous, composing it with mollified
$u$ and $v$ would not yield a smooth coefficient; we therefore mollify $V$ itself. Let
$\rho\in C^\infty_0(\R^{n+1})$, $\rho\ge0$, $\int\rho=1$, and $\rho_k(z):=k^{n+1}\rho(kz)$. If $\tilde V$ denotes the
extension of $V$ by zero outside $(\tau,T)\times\Om$, set
\begin{equation}\label{mollified}
  V_k:=\big(\tilde V*\rho_k\big)\big|_{[\tau,T]\times\overline\Om} .
\end{equation}
The mollification is therefore taken in space \emph{and} in time. For each $k$, consider the backward approximate
problem
\begin{equation}\label{dual-k}
  \Dtt\vp_k-\Dx\vp_k+(V_k-m)\vp_k=F\ \text{ in }(\tau,T)\times\Om,\quad
  \vp_k|_{\partial\Om}=0,\quad \vp_k(T)=\Dt\vp_k(T)=0 .
\end{equation}

Since $\tilde V\in L^\mu(\R^{n+1})$, Young's inequality gives, for every multi-index $\beta$,
$|\partial^\beta_{t,x}(\tilde V*\rho_k)|\le\|\partial^\beta\rho_k\|_{L^{\mu'}(\R^{n+1})}
  \|\tilde V\|_{L^{\mu}(\R^{n+1})}$, so $V_k$ and all its
derivatives are bounded on $[\tau,T]\times\overline\Om$ for each fixed $k$. Moreover, Minkowski's inequality yields
\[
  \|(\tilde V*\rho_k)(t)\|_{L^\mu(\R^n)}
  \le\iint_{\R^{n+1}}\rho_k(s,y)\,\|\tilde V(t-s)\|_{L^\mu(\R^n)}\,\td s\,\td y,
\]
and therefore, by \eqref{V-R},
\begin{equation}\label{Vk-R}
  \sup_k\|V_k\|_{L^\infty(\tau,T;L^{\mu}(\Om))}\le\|V\|_{L^\infty(\tau,T;L^{\mu}(\Om))}
  \le C_1\(1+2R^{p-1}\).
\end{equation}

\begin{lemma}\label{lem-Vk}
  For every $\vartheta\in[1,\infty)$,
  $\displaystyle\lim_{k\to\infty}\|V-V_{k}\|_{L^{\vartheta}(\tau,T;L^{\mu}(\Om))}=0$.
\end{lemma}

\begin{proof}
  Since $p>1$, we have $\mu<\infty$, and $\tilde V\in L^{\vartheta}(\R;L^{\mu}(\R^n))$ because
  $\tilde V$ is bounded in $L^\mu(\R^n)$ uniformly in $t$ and vanishes for $t\notin(\tau,T)$. As
  $\vartheta,\mu<\infty$, compactly supported continuous functions are dense in this mixed
  space, so translations $z\mapsto\tilde V(\cdot-z)$ are continuous in it. Minkowski's
  inequality gives
  \[
    \|\tilde V*\rho_k-\tilde V\|_{L^{\vartheta}(\R;L^{\mu}(\R^n))}
    \le\int_{\R^{n+1}}\rho_k(z)\,\|\tilde V(\cdot-z)-\tilde V\|_{L^{\vartheta}(\R;L^{\mu}(\R^n))}\,\td z
    \longrightarrow0,
  \]
  and restricting to $(\tau,T)\times\Om$, where $\tilde V=V$, proves the claim.
\end{proof}

We now solve the regularised problem \eqref{dual-k} with a uniform Strichartz regularity. Choose
\begin{equation}
  (q,r)=
  \begin{cases}
    \Big(\dfrac{7}{2},\ 14\Big),         & n=3,   \\[3mm]
    \Big(n-1,\ \dfrac{2(n-1)}{n-3}\Big), & n\ge4,
  \end{cases}
\end{equation}
so that $(q,r,1)$ is $\Om$-admissible; Appendix \ref{app-opt} derives this pair.
\begin{prop}\label{prop-dualbound}
  For each fixed $k$, problem \eqref{dual-k} has a unique solution $\vp_k$ such that
  \begin{equation}\label{dualreg}
    \vp_k\in C^2\([\tau,T];L^2(\Om)\)\cap C^1\([\tau,T];H^1_0(\Om)\)
    \cap C\([\tau,T];H^2(\Om)\cap H^1_0(\Om)\).
  \end{equation}
  Moreover there exists $\ell=\ell(R,\ga,n)\in(0,1]$, which depends on the two solutions only through $R$, such
  that whenever $T-\tau\le\ell$,
  \begin{equation}\label{dualbound}
    \sup_k\|\vp_k\|_{L^{q}(\tau,T;L^{r}(\Om))}\le 2C_0\|F\|_{L^1(\tau,T;L^2(\Om))}<\infty.
  \end{equation}
\end{prop}

\begin{proof}
  Reversing time, \eqref{dual-k} is a linear wave equation with Dirichlet boundary condition,
  vanishing initial data, and a coefficient $V_k-m$ which is smooth and bounded on
  $[\tau,T]\times\overline\Om$, driven by $F\in C^\infty_0((\tau,T)\times\Om)$. Since the data
  vanish identically and $F(t,\cdot)\in C^\infty_0(\Om)$, all compatibility conditions at the
  initial time are satisfied, so the classical energy method in the scale $H^2\cap H^1_0$
  applies and gives \eqref{dualreg}.

  Treating $(V_k-m)\vp_k$ as a source and applying the Strichartz estimate \eqref{BSS-inhom} on $(\tau,T)$ to
  \eqref{dual-k} gives
  \begin{equation}\label{dual-str}
    \|\vp_k\|_{L^{q}(\tau,T;L^{r}(\Om))}
    \le C_0\Big(\|F\|_{L^{1}(\tau,T;L^2(\Om))}+\|(V_k-m)\vp_k\|_{L^{1}(\tau,T;L^2(\Om))}\Big).
  \end{equation}
  Since $\frac1\mu+\frac1r=\frac12$, by H\"older inequality and \eqref{Vk-R},
  \begin{equation}\label{holder-concrete}
    \begin{aligned}
      \|(V_k-m)\vp_k\|_{L^{1}(\tau,T;L^2(\Om))}
       & \le(T-\tau)^{1/{q'}}\|V_k-m\|_{L^{\infty}(\tau,T;L^{\mu}(\Om))}\|\vp_k\|_{L^{q}(\tau,T;L^{r}(\Om))} \\
       & \le (T-\tau)^{1/{q'}}\bigl(C_1(1+2R^{p-1})+m|\Om|^{1/\mu}\bigr)
      \|\vp_k\|_{L^{q}(\tau,T;L^{r}(\Om))},
    \end{aligned}
  \end{equation}
  where $q'$ is the H\"older conjugate of $q$. Choose $\ell\in(0,1]$ so that
  \begin{equation}\label{ell-choice}
    \ell^{1/q'}C_0\bigl(C_1(1+2R^{p-1})+m|\Om|^{1/\mu}\bigr)\le\frac12.
  \end{equation}
  If $T-\tau\le\ell$, the last term in \eqref{dual-str} is absorbed into the left-hand side, which yields
  \eqref{dualbound}. Note that $\ell$ involves only $R,\ga,\Om,n$, so in particular it is independent of $k$ and of $F$.
\end{proof}

\begin{remark}\label{rem-nosupport}
  Two statements that one might expect are neither needed nor available. No bound of $\vp_k$ in
  $L^\infty((\tau,T)\times\Om)$ is asserted: the $L^qL^r$ bound is exactly what the duality argument requires. Nor is
  $\vp_k$ compactly supported in $\Om$, since on a bounded domain with the Dirichlet condition the waves reach and are
  reflected off $\partial\Om$, so that finite speed of propagation yields no such conclusion. The admissibility of
  $\vp_k$ as a test function is settled by the $H^{-1}$-valued form of \eqref{delta-eq}.
\end{remark}

\begin{proof}[Proof of Theorem \ref{thm-unique}]
  Let $\ell$ be as in Proposition \ref{prop-dualbound} and suppose first $T-\tau\le\ell$. Fix
  $F\in C^\infty_0((\tau,T)\times\Om)$ and let $\vp_k$ solve \eqref{dual-k}.

  Equation \eqref{delta-eq} holds in $L^{\infty}(\tau,T;H^{-1}(\Om))$, so it may be paired with $\vp_k(t)$ and
  integrated:
  \begin{equation}
    \int_\tau^T\Big(\big\la\de_{tt},\vp_k\big\ra+\big\la\Nx\de,\Nx\vp_k\big\ra
    +\big\la(V-m)\de,\vp_k\big\ra\Big)\,\td t=0,
  \end{equation}
  all three terms being finite, the last by \eqref{Vdelta} together with
  $\vp_k\in C_tH^1_0\subset C_tL^{2^*}$ and $\frac{p}{2^*}+\frac1{2^*}\le1$, which is again the condition
  $p\le p_c$. Since $\Dt\de\in L^{\infty}(\tau,T;L^2)$ and $\Dtt\de\in L^{\infty}(\tau,T;H^{-1})$, two
  integrations by parts in time give
  \[
    \int_\tau^T\big\la\de_{tt},\vp_k\big\ra\,\td t
    =\Big[\big\la\Dt\de,\vp_k\big\ra-\big\la\de,\Dt\vp_k\big\ra\Big]_{\tau}^{T}
    +\int_\tau^T\big\la\de,\Dtt\vp_k\big\ra\,\td t
    =\int_\tau^T\big\la\de,\Dtt\vp_k\big\ra\,\td t,
  \]
  because $\de(\tau)=\Dt\de(\tau)=0$ and $\vp_k(T)=\Dt\vp_k(T)=0$. Finally,
  $\vp_k\in L^{\infty}(\tau,T;H^2\cap H^1_0)$ and $\de\in L^{\infty}(\tau,T;H^1_0)$ allow the transfer
  $\int\la\Nx\de,\Nx\vp_k\ra=\int\la\de,-\Dx\vp_k\ra$. Combining these identities gives
  \[
    0=\int_\tau^T\!\!\int_\Om\de\Big(\Dtt\vp_k-\Dx\vp_k+(V-m)\vp_k\Big)\td x\,\td t
    =\int_\tau^T\!\!\int_\Om\de\Big(F+(V-V_k)\vp_k\Big)\td x\,\td t,
  \]
  where \eqref{dual-k} was used. Hence
  \begin{equation}\label{almost}
    \left|\int_\tau^T\!\!\int_\Om\de F\,\td x\,\td t\right|
    \le\|\de\|_{L^\infty(\tau,T;L^{2^*})}\,
    \big\|(V-V_k)\vp_k\big\|_{L^{1}(\tau,T;L^{(2^*)'})}.
  \end{equation}
  Applying H\"older's inequality with $\frac1\mu+\frac1r=\frac12$ in space and
  $\frac{1}{q'}+\frac1q=1$ in time together with
  \eqref{dualbound} and $T-\tau\le1$, we obtain
  \[
    \begin{aligned}
      \big\|(V-V_k)\vp_k\big\|_{L^{1}(\tau,T;L^2(\Om))}
       & \le\|V-V_k\|_{L^{q'}(\tau,T;L^{\mu}(\Om))}\,\|\vp_k\|_{L^{q}(\tau,T;L^{r}(\Om))} \\
       & \le2C_0\|F\|_{L^1(\tau,T;L^{2}(\Om))}\,\|V-V_k\|_{L^{q'}(\tau,T;L^{\mu}(\Om))}.
    \end{aligned}
  \]
  Letting $k\to\infty$ and applying Lemma \ref{lem-Vk} with $\vartheta=q'$, using
  $L^2(\Om)\hookrightarrow L^{(2^*)'}(\Om)$, shows that the right-hand side of \eqref{almost} tends to $0$, whereas
  the left-hand side does not depend on $k$. Therefore
  \[
    \int_\tau^T\!\!\int_\Om\de F\,\td x\,\td t=0
    \qquad\text{for every }F\in C^\infty_0\((\tau,T)\times\Om\),
  \]
  whence $\de=0$ almost everywhere, that is, $u=v$ on $[\tau,T]$.

  For general $T>\tau$, note that $\ell$ in \eqref{ell-choice} depends only on the bound $R$, which is finite on all
  of $[\tau,T]$. Partition $[\tau,T]$ into finitely many intervals of length at most $\ell$. On the first of them the
  argument above gives $u=v$, hence $w(t_1)=\Dt w(t_1)=0$, and it may then be repeated verbatim with $\tau$ replaced
  by $t_1$, the same $\ell$ being admissible. After finitely many steps $u=v$ on $[\tau,T]$; since $T$ was
  arbitrary, $u=v$ on $[\tau,\infty)$.
\end{proof}

A weak solution on $[s,t]$ satisfies $u\in L^\infty(s,t;H^1_0)$, $\Dt u\in L^\infty(s,t;L^2)$ and, by the
equation, $\Dtt u\in L^\infty(s,t;H^{-1})$. Hence $u\in C_w([s,t];H^1_0(\Om))$ and
$\Dt u\in C_w([s,t];L^2(\Om))$, so that $\xi_u(\varsigma)\in\H$ is defined for every $\varsigma\in[s,t]$.
We say that $u$ is a weak solution of \eqref{eq1} on $[s,t]$ if it satisfies Definition
\ref{def-weak} with $(\tau,T)$ replaced by $(s,t)$, without prescribing its data at $s$.

\begin{cor}[Backward uniqueness]\label{cor-back-unique}
  Assume \eqref{F1}. If $u,v$ are weak solutions of \eqref{eq1} on $[s,t]$ with
  $\xi_u(t)=\xi_v(t)$, then $u=v$ on $[s,t]$.
\end{cor}

\begin{proof}
  The equation \eqref{delta-eq}, with $\tau$ replaced by $s$, contains no first-order time
  derivative and is therefore invariant under time reversal $\varsigma\mapsto s+t-\varsigma$. Under this
  transformation, the potential $V(\varsigma,x)-m$ becomes $V(s+t-\varsigma,x)-m$, which still satisfies
  the required bounds: $V\in L^\infty(\R;L^\mu(\Om))$ by \eqref{V-R} is invariant under
  time shifts, and $m$ is a constant.

  Now $\de(t)=\Dt\de(t)=0$, and the dual problems \eqref{dual-k} are posed forward in (reversed)
  time with $\vp_k(s)=\Dt\vp_k(s)=0$. Proposition \ref{prop-dualbound}, whose proof uses only the
  Strichartz estimate on $(s,t)$ and the $L^\infty L^\mu$ bound on $V$, applies without change.
  In the integration by parts of the proof of Theorem \ref{thm-unique} the boundary terms now
  vanish at $\varsigma=t$ because of the data of $\de$, and at $\varsigma=s$ because of the data of $\vp_k$.
  The rest of the argument, including the iteration over intervals of length $\ell$ (now starting
  from $t$ and moving backwards), is unchanged.
\end{proof}

\subsection{Existence and the extra space-time regularity}\label{ssec-exist}
\begin{thm}\label{thm-exist}
  Assume \eqref{F1}, \eqref{F2} and \eqref{G1}. Then for every $\xi_\tau\in\H$ problem \eqref{eq1} has a unique weak
  solution $u$, and there are $\al>0$ and a monotone function $Q$, depending on the structural parameters $\ga,\lm,\lm_1,p,n,\Om$, such that
  \begin{equation}\label{extra-reg}
    \|\xi_u(t)\|_{\H}+\big\|u-\Theta\big\|_{L^{q}(t,t+1;L^{r}(\Om))}
    \le Q\(\|\xi_\tau\|_{\H}+\|g\|_{H^1_b}\)e^{-\al(t-\tau)}+Q\(\|g\|_{H^1_b}\),\quad t\ge\tau,
  \end{equation}
  where $\Theta$ is the unique solution of \eqref{eq1} with vanishing initial data.
\end{thm}

\begin{proof}
  \emph{Step 1: approximation.} Let $P_N$ be the orthogonal projection onto $\mathrm{span}\{e_1,\dots,e_N\}$. Consider
  \begin{equation}\label{galerkin}
    \Dtt u_N+\ga\Dt u_N-\Dx u_N+P_Nf(u_N)=P_Ng,\qquad
    \xi_{u_N}(\tau)=\(P_Nu_0,P_Nu_1\),
  \end{equation}
  an ordinary differential system with locally Lipschitz right-hand side, hence locally
  solvable.

  Set
  \[
    \begin{aligned}
      E(u)        & :=\frac12\|\Dt u\|^2+\frac12\|\Nx u\|^2+\int_\Om F(u)\,\td x-\big\la g,u\big\ra, \\
      \Phi_\al(u) & :=E(u)+\al\big\la u,\Dt u\big\ra+\frac{\al\ga}{2}\|u\|^2.
    \end{aligned}
  \]

  Testing the equation \eqref{galerkin} with $\Dt u_N+\al u_N$ gives
  \begin{multline}\label{energy-ineq}
    \frac{\td}{\td t}\Phi_\al(u_N)+\al\Phi_\al(u_N)
    =-\Big(\ga-\frac{3\al}{2}\Big)\|\Dt u_N\|^2-\frac{\al}{2}\|\Nx u_N\|^2
    -\al\int_\Om\big[f(u_N)u_N-F(u_N)\big]\td x\\
    -\big\la\Dt g,u_N\big\ra+\al^2\la u_N,\Dt u_N\ra+\frac{\al^2\ga}{2}\|u_N\|^2 .
  \end{multline}
  By the third inequality of \eqref{F2} and $\|u_N\|^2\le\lm_1^{-1}\|\Nx u_N\|^2$,
  \[
    -\al\int_\Om\big[f(u_N)u_N-F(u_N)\big]\td x
    \le\al\frac{\lm}{2\lm_1}\|\Nx u_N\|^2+\al C_f|\Om|.
  \]
  Moreover
  $$-\la\Dt g,u_N\ra\le\frac C\al\|\Dt g\|^2_{H^{-1}}+\frac{\lm_1-\lm}{4\lm_1}\al\|\Nx u_N\|^2,$$
  and the last two terms
  of \eqref{energy-ineq} are bounded by
  $\frac{\al^2}{2}\|\Dt u_N\|^2+\frac{\al^2(1+\ga)}{2\lm_1}\|\Nx u_N\|^2$. Choosing $\al>0$ with
  $$\ga-\frac{3\al}{2}-\frac{\al^2}{2}\ge0 \text{\quad and\quad} \frac{\al^2(1+\ga)}{2\lm_1}\le\frac{\lm_1-\lm}{4\lm_1}\al.$$
  Then
  \begin{equation}\label{energy-final}
    \frac{\td}{\td t}\Phi_\al(u_N)+\al\Phi_\al(u_N)
    \le \frac C\al\|\Dt g\|^2_{H^{-1}}+\al C_f|\Om| .
  \end{equation}
  Conversely, by the second inequality of \eqref{F2},
  $$\int F(u_N)\ge-\frac{\lm}{2\lm_1}\|\Nx u_N\|^2-C_f|\Om|,$$
  while
  $$-\la g,u_N\ra\ge-\frac{\lm_1-\lm}{4\lm_1}\|\Nx u_N\|^2-C\|g\|^2_{H^{-1}}$$
  and
  $$\al\la u_N,\Dt u_N\ra\ge-\frac{\al}{2\lm_1}\|\Dt u_N\|^2-\frac{\al}{2}\|\Nx u_N\|^2.$$
  Hence for $\al < \min\{\frac12, \frac{\lm_1-\lm}{4}\}$,
  \begin{equation}\label{Phi-equiv}
    \frac\al2\|\xi_{u_N}\|^2_{\H}-C\(\|g\|^2_{H^{-1}}+1\)
    \le\Phi_\al(u_N)\le Q\(\|\xi_{u_N}\|_{\H}\)+C\|g\|^2_{H^{-1}},
  \end{equation}
  the upper bound using $|F(s)|\le C(1+|s|^{p+1})$ from \eqref{fbounds} together with
  $H^1_0(\Om)\hookrightarrow L^{p+1}(\Om)$, legitimate since $p+1<2^*$.

  Note that
  $$\| g\|_{H^{-1}}\le \|g\|_{H^1_b(\R;H^{-1})}.$$
  Applying Gronwall's Lemma to \eqref{energy-final} implies
  \begin{equation}\label{apr-energy}
    \sup_{N}\ \sup_{t\ge\tau}\|\xi_{u_N}(t)\|_{\H}\le Q\(\|\xi_\tau\|_{\H}+\|g\|_{H^1_b(\R;H^{-1})}\)e^{-\al(t-\tau)}+Q\(\|g\|_{H^1_b(\R;H^{-1})}\).
  \end{equation}
  Set
  \begin{equation}
    R_0= Q\(\|\xi_\tau\|_{\H}+\|g\|_{H^1_b(\R;H^{-1})}\)+Q\(\|g\|_{H^1_b(\R;H^{-1})}\),
  \end{equation}
  then
  \begin{equation}\label{R-def}
    \sup_{N}\ \sup_{t\ge\tau}\|\xi_{u_N}(t)\|_{\H}\le R_0.
  \end{equation}
  Let $\Theta_N$ solve \eqref{galerkin} with vanishing data; then also
  $$\sup_{t\ge\tau}\|\xi_{\Theta_N}(t)\|_{\H}\le Q\(\|g\|_{H^1_b(\R;H^{-1})}\).$$

  \emph{Step 2: Strichartz on short intervals.} Put $v_N:=u_N-\Theta_N$. Subtracting the two Galerkin
  systems, the term $P_Ng$ cancels \emph{identically} and
  \begin{equation}\label{vN-eq}
    \Dtt v_N+\ga\Dt v_N-\Dx v_N=P_N\big[f(\Theta_N)-f(u_N)\big]=-P_N\(W_Nv_N\),
    \quad \xi_{v_N}(\tau)=\(P_Nu_0,P_Nu_1\),
  \end{equation}
  where $W_N:=\int_0^1f'\(su_N+(1-s)\Theta_N\)\td s$ satisfies, by \eqref{V-R} and \eqref{R-def},
  \begin{equation}\label{WN}
    \|W_N\|_{L^\infty(\tau,\infty;L^{\mu}(\Om))}\le C_1\(1+2R_0^{p-1}\)=:\kappa_0 .
  \end{equation}

  Fix any $t\ge\tau$. By \eqref{apr-energy}, $\|\xi_{v_N}(t)\|_{\H}\le \|\xi_{u_N}(t)\|_{\H}+\|\xi_{\Theta_N}(t)\|_{\H}$ satisfies
  \begin{equation}\label{vN-energy-t}
    \|\xi_{v_N}(t)\|_{\H}\le Q\(\|\xi_\tau\|_{\H}+\|g\|_{H^1_b}\)e^{-\al(t-\tau)}+Q\(\|g\|_{H^1_b}\).
  \end{equation}
  Let $I=[t,t+\ell]$, $\ell\in(0,1]$. By Strichartz estimate \eqref{lem-damped} applied to \eqref{vN-eq} on $I$,
  $\|P_N\(W_Nv_N\)\|\le\|W_Nv_N\|$ and \eqref{holder-concrete},
  \[
    \|v_N\|_{L^{q}(I;L^{r}(\Om))}
    \le C_0\Big(\|\xi_{v_N}(t)\|_{\H}+\ell^{1/{q'}}\kappa_{0}\|v_N\|_{L^{q}(I;L^{r}(\Om))}\Big).
  \]
  Fix $\ell=\ell(\kappa_0)\in(0,1]$ with $C_0\ell^{1/{q'}}\kappa_{0}\le\frac12$; the left-hand side is
  finite because $v_N$ is a finite sum of smooth modes, so by \eqref{vN-energy-t},
  \begin{equation}\label{vN-str}
    \|v_N\|_{L^{q}(t,t+\ell;L^{r}(\Om))}\le 2C_0\|\xi_{v_N}(t)\|_{\H}
    \le 2C_0\left(Q\(\|\xi_\tau\|_{\H}+\|g\|_{H^1_b}\)e^{-\al(t-\tau)}+Q\(\|g\|_{H^1_b}\)\right).
  \end{equation}
  Covering $[t,t+1]$ by $\lceil1/\ell\rceil$ such intervals and summing the $q$-th powers,
  \begin{equation}\label{vN-sum}
    \|v_N\|_{L^{q}(t,t+1;L^{r}(\Om))}\le\lceil1/\ell\rceil^{1/q}\,2C_0\left(Q\(\|\xi_\tau\|_{\H}+\|g\|_{H^1_b}\)e^{-\al(t-\tau)}+Q\(\|g\|_{H^1_b}\)\right),
  \end{equation}
  uniformly in $N$.

  \emph{Step 3: passage to the limit.} By \eqref{apr-energy} the sequences $u_N,\Dt u_N$ are bounded
  in $L^\infty(\tau,T;H^1_0)$, $L^\infty(\tau,T;L^2)$, and from \eqref{galerkin} $\Dtt u_N$ is
  bounded in $L^\infty(\tau,T;H^{-1})$, since
  $\|P_Nf(u_N)\|_{H^{-1}}\le C\|f(u_N)\|_{L^{(2^*)'}}\le C(1+R_0^{p})$ by
  $\frac{2^*}{p}>(2^*)'$. By the Aubin--Lions--Simon theorem \cite{simon87} a subsequence
  satisfies
  \begin{equation}\label{strongconv}
    u_N\to u\ \text{ in }C\([\tau,T];H^{1-\de}(\Om)\)\ \forall\de\in(0,1],\qquad
    \xi_{u_N}\overset{*}{\wto}\xi_u\ \text{ in }L^\infty(\tau,T;\H).
  \end{equation}
  Since $H^{1-\de}(\Om)\hookrightarrow L^{\eta}(\Om)$ for some $\eta<2^*$ with $\eta\to2^*$ as
  $\de\to0$, and since $p<p_c$ leaves a margin in the H\"older count, the bound
  $|f(u_N)-f(u)|\le C\(1+|u_N|^{p-1}+|u|^{p-1}\)|u_N-u|$ together with the uniform
  $L^{2^*}$-bound gives $f(u_N)\to f(u)$ in $L^{\infty}(\tau,T;L^{\vartheta}(\Om))$ for some
  $\vartheta>(2^*)'$. Hence $P_Nf(u_N)\to f(u)$ in $L^\infty(\tau,T;H^{-1})$ and $u$ is a weak
  solution. By Theorem \ref{thm-unique} the
  limit is independent of the subsequence, so $u$ is \emph{the} unique weak solution. The same
  argument gives $\Theta_N\to\Theta$, where $\Theta$ is the unique solution of \eqref{eq1} with vanishing initial data
  \begin{align}
    \begin{cases}
      \Theta_{tt}+\ga \Theta_{t}-\Dx \Theta+f(\Theta)=g(t), ~(t,x)\in (\tau, +\infty)\times\Omega, \\
      \Theta(\tau,x)=0,~~ \Dt\Theta(\tau,x)=0,  x\in \Omega,                                       \\
      \Theta(t,x)=0,~ (t,x)\in (\tau, +\infty)\times \partial\Omega.
    \end{cases}
  \end{align}
  Hence $v_N\to u-\Theta$ in $L^{\infty}(\tau,T;H^{1-\de})$.

  Finally $L^{q}(t,t+1;L^{r}(\Om))$ is reflexive, so by \eqref{vN-sum} a further subsequence of
  $v_N$ converges weakly in it, necessarily to $u-\Theta$, and the norm is weakly lower
  semicontinuous:
  \[
    \begin{aligned}
      \|u-\Theta\|_{L^{q}(t,t+1;L^{r}(\Om))} & \le\liminf_N\|v_N\|_{L^{q}(t,t+1;L^{r}(\Om))}                                                                           \\
                                             & \le \lceil1/\ell\rceil^{1/q}\,2C_0\left(Q\(\|\xi_\tau\|_{\H}+\|g\|_{H^1_b}\)e^{-\al(t-\tau)}+Q\(\|g\|_{H^1_b}\)\right).
    \end{aligned}
  \]
  With \eqref{apr-energy} this gives \eqref{extra-reg}.
\end{proof}

\begin{remark}\label{rem-honest}
  The second term in \eqref{extra-reg} involves $u-\Theta$ and not $u$, and this is not an artefact of the
  proof. Under \eqref{G1} the function $\Theta$ is itself known only to lie in the energy class, so \eqref{extra-reg}
  yields \emph{no} information about $\|u\|_{L^{q}(L^{r})}$. What \eqref{extra-reg} does provide, and what is used
  below, is control of the difference of two solutions: applying it to two sets of data and subtracting eliminates
  $\Theta$ altogether, which is Theorem \ref{thm-lip}.
\end{remark}

\subsection{Lipschitz continuous dependence}\label{ssec-lip}

\begin{thm}\label{thm-lip}
  Assume \eqref{F1}, \eqref{F2} and \eqref{G1}. Let $u^1,u^2$ be the weak solutions with data
  $\xi^1_\tau,\xi^2_\tau\in\H$, and put $\bar u:=u^1-u^2$. Then there exist $K,C>0$, depending on the structural
  parameters and on $R_0:=\max_iQ\bigl(\|\xi^i_\tau\|_{\H}+\|g\|_{H^1_b}\bigr)$, such that
  \begin{equation}\label{lip}
    \|\xi_{\bar u}(t)\|_{\H}+\|\bar u\|_{L^{q}(\tau,t;L^{r}(\Om))}
    \le Ce^{K(t-\tau)}\|\xi_{\bar u}(\tau)\|_{\H},\qquad t\ge\tau.
  \end{equation}
\end{thm}

\begin{proof}
  As in \eqref{w-eq}, the difference $\bar u$ satisfies
  $\Dtt\bar u+\ga\Dt\bar u-\Dx\bar u=-V\bar u$ with $\|V\|_{L^\infty_tL^{\mu}}\le\kappa_0$, the force having
  cancelled. Let $\ell\in(0,1]$ satisfy $C_0\ell^{1/q'}\kappa_0\le\frac12$, and let $I=[s,s+\ell]$. Exactly as in
  Step 2 above,
  \begin{equation}\label{lip-str}
    \|\bar u\|_{L^{q}(I;L^{r}(\Om))}\le 2C_0\|\xi_{\bar u}(s)\|_{\H},
  \end{equation}
  the left-hand side being finite by Theorem \ref{thm-exist}, since
  $\bar u=(u^1-\Theta)-(u^2-\Theta)$ is a difference of two functions each controlled by \eqref{extra-reg}.
  Consequently, by \eqref{holder-concrete},
  \begin{equation}\label{lip-source}
    \|V\bar u\|_{L^1(I;L^2)}\le\ell^{1/q'}\kappa_0\|\bar u\|_{L^{q}(I;L^{r}(\Om))}
    \le\frac1{2C_0}\cdot2C_0\|\xi_{\bar u}(s)\|_{\H}=\|\xi_{\bar u}(s)\|_{\H}.
  \end{equation}
  Testing the equation with $\Dt\bar u$ on the other hand gives
  \[
    \frac{\td}{\td t}\Big(\frac12\|\Dt\bar u\|^2+\frac12\|\Nx\bar u\|^2\Big)
    +\ga\|\Dt\bar u\|^2=-\big\la V\bar u,\Dt\bar u\big\ra\le\|V\bar u\|\,\|\Dt\bar u\|,
  \]
  so $\frac{\td}{\td t}\|\xi_{\bar u}\|_{\H}\le\|V\bar u\|$ and, by \eqref{lip-source},
  \begin{equation}\label{lip-energy}
    \|\xi_{\bar u}(t)\|_{\H}\le\|\xi_{\bar u}(s)\|_{\H}+\|V\bar u\|_{L^1(I;L^2)}
    \le2\|\xi_{\bar u}(s)\|_{\H},\qquad t\in I.
  \end{equation}
  Since $\ell$ depends only on $R_0$, which by Theorem \ref{thm-exist} bounds both solutions uniformly on
  $[\tau,\infty)$, the estimates \eqref{lip-str} and \eqref{lip-energy} iterate over the
  $\lceil(t-\tau)/\ell\rceil$ successive intervals, giving
  $\|\xi_{\bar u}(t)\|_{\H}\le2^{(t-\tau)/\ell+1}\|\xi_{\bar u}(\tau)\|_{\H}$; summing the $q$-th powers of
  \eqref{lip-str} bounds $\|\bar u\|_{L^{q}(\tau,t;L^{r}(\Om))}$. This is \eqref{lip} with
  $K=2\ln2/\ell$.
\end{proof}

\begin{remark}\label{rem-two-uniq}
  Theorems \ref{thm-unique} and \ref{thm-lip} are logically distinct, and neither implies the other.
  Theorem \ref{thm-unique} is \emph{unconditional}: it compares two arbitrary energy-class solutions and uses no
  space-time integrability. This is precisely what makes $\Theta$ well defined, and hence what renders the present
  section non-circular. Theorem \ref{thm-lip} is quantitative, but it is stated for the solutions produced by
  Theorem \ref{thm-exist}, whose Strichartz regularity relative to $\Theta$ is what makes the left-hand side of
  \eqref{lip-str} finite in the first place. Conflating the two statements would make the argument circular.
\end{remark}


\subsection{Exponentially weighted energy identity}\label{ssec-expenergy}

\begin{lemma}\label{lem-expenergy}
  Assume \eqref{F1}, \eqref{F2} and \eqref{G1}, and let $u$ be a weak solution of \eqref{eq1} on $[s_0,t_0]$.
  Define
  \begin{align}
    \mathscr E(t) & :=e^{\ga t}\Bigl(\frac12\|\Dt u(t)+\tfrac\ga2 u(t)\|^2+\frac12\|\Nx u(t)\|^2
    -\frac{\ga^2}{8}\|u(t)\|^2+\!\int_\Om F\bigl(u(t)\bigr)\,\td x
    -\la g(t),u(t)\ra\Bigr),\label{expE}                                                         \\
    R(t)          & :=\ga e^{\ga t}\!\int_\Om\Bigl[\frac12f(u)u-F(u)\Bigr](t)\,\td x
    +e^{\ga t}\Bigl\la\Dt g(t)+\tfrac\ga2g(t),u(t)\Bigr\ra .\label{expR}
  \end{align}
  Then $R\in L^1(s_0,t_0)$, and for all $s_0\le s\le t\le t_0$,
  \begin{equation}\label{expid}
    \mathscr E(t)+\int_s^tR(\varsigma)\,\td\varsigma =\mathscr E(s).
  \end{equation}
\end{lemma}

The point of the weight is that \eqref{expid} contains no dissipation integral $\int\|\Dt u\|^2$. Such an
integral is only lower semicontinuous under weak convergence, and would therefore prevent passage to the limit in
one of the two time directions. After the substitution, the damping appears only through the endpoint quantity
$\mathscr E$ and through $R$, which does not involve $\Dt u$.

\begin{proof}
  All terms of $\mathscr E$ and $R$ are finite. For $F(u)$ and $f(u)u$ this follows from
  \eqref{fbounds} and $p+1<2^*$; for the pairings, from $g,\Dt g\in L^2_{loc}(H^{-1})$ and
  $u\in L^\infty(H^1_0)$. In particular $R\in L^1(s_0,t_0)$. Fix $s_0\le s<t\le t_0$.

  \emph{Step 1: the identity for Galerkin solutions.} Let $u_N$ be the solution of the
  system \eqref{galerkin} on $[s,t]$, with data prescribed at $s$ or at $t$, and put
  $\de_N:=e^{\ga\varsigma/2}u_N$. A direct computation gives
  \[
    \partial_{\varsigma\varsigma}\de_N-\Dx\de_N-\frac{\ga^2}{4}\de_N=e^{\ga\varsigma/2}\(P_Ng-P_Nf(u_N)\).
  \]
  Take the $L^2$ inner product with $\partial_\varsigma\de_N$. Using $\partial_\varsigma\de_N=e^{\ga\varsigma/2}\(\partial_\varsigma u_N+\frac\ga2u_N\)$,
  the chain rule for $F$, and the absolute continuity of $\varsigma\mapsto\la e^{\ga\varsigma/2}g,\de_N\ra$
  (valid since $e^{\ga\varsigma/2}g\in H^1(s,t;H^{-1})$ and $\de_N\in C^1([s,t];H^1_0)$), we get
  \[
    \big\la e^{\ga\varsigma/2}f(u_N),\partial_\varsigma\de_N\big\ra
    =\frac{\td}{\td\varsigma}\Big[e^{\ga\varsigma}\!\int_\Om F(u_N)\Big]
    -\ga e^{\ga\varsigma}\!\int_\Om F(u_N)+\frac\ga2e^{\ga\varsigma}\la f(u_N),u_N\ra,
  \]
  \[
    \la e^{\ga\varsigma/2}g,\partial_\varsigma\de_N\ra=\frac{\td}{\td\varsigma}\la e^{\ga\varsigma/2}g,\de_N\ra-\la\partial_\varsigma (e^{\ga\varsigma/2}g),\de_N\ra,\qquad
    \la\partial_\varsigma (e^{\ga\varsigma/2}g),\de_N\ra=e^{\ga\varsigma}\Big\la\partial_\varsigma g+\frac\ga2g,u_N\Big\ra .
  \]
  Hence $\mathscr E_N$ and $R_N$, defined by \eqref{expE}--\eqref{expR} with $u_N$ in place of
  $u$, satisfy \emph{exactly}
  \begin{equation}\label{gal-exp}
    \mathscr E_N(t)+\int_s^tR_N(\varsigma)\,\td\varsigma=\mathscr E_N(s).
  \end{equation}

  \emph{Step 2: Galerkin solutions from either endpoint.} Let $\varsigma_0\in\{s,t\}$ and let $u_N$ be
  the Galerkin solution with data $P_N\xi_u(\varsigma_0)$ at $\varsigma_0$. We show that it exists on all of
  $[s,t]$ with a bound independent of $N$.

  If $\varsigma_0=s$, this is \eqref{extra-reg} with $\tau=s$.

  If $\varsigma_0=t$, we solve backwards and use the unweighted energy
  $E(u_N)=\frac12\|\Dt u_N\|^2+\frac12\|\Nx u_N\|^2+\int F(u_N)-\la g,u_N\ra$. Multiplying the
  Galerkin equation by $\Dt u_N$ and integrating gives, exactly as for Theorem \ref{thm-exist},
  \[
    E(u_N(s))=E(u_N(t))+\int_s^t\Big(\ga\|\Dt u_N\|^2+\la\Dt g,u_N\ra\Big)\td\varsigma .
  \]
  By \eqref{F2} and the estimate $|\la g,u_N\ra|\le\frac{\al}{\lm_1}\|\Nx u_N\|^2+C\|g\|^2_{H^{-1}}$ with $\al = \min\{\frac12, \frac{\lm_1-\lm}{4}\}$,
  $$E(u_N)\ge \al\|\xi_{u_N}\|^2_{\H}-C_g,$$
  with
  $$C_g:=C\(\|g\|^2_{L^\infty(s,t;H^{-1})}+1\).$$
  Hence $\Psi_N:=E(u_N)+C_g\ge\al\|\xi_{u_N}\|^2_{\H}\ge 0$, and
  \[
    \Psi_N(s)\le\Psi_N(t)
    +\int_s^t\Big(C\Psi_N+\|\Dt g\|^2_{H^{-1}}\Big)\td\varsigma .
  \]
  Gronwall's inequality backwards in time gives
  \begin{equation}\label{back-gron}
    \|\xi_{u_N}(\varsigma)\|^2_{\H}\le4\Big(Q\(\|\xi_u(t)\|_{\H}\)+C_g
    +\|\Dt g\|^2_{L^2(s,t;H^{-1})}\Big)e^{C(t-\varsigma)},\qquad \varsigma\in[s,t],
  \end{equation}
  using $\|P_N\xi_u(t)\|_{\H}\le\|\xi_u(t)\|_{\H}$. In particular the backward solution does not
  blow up and exists on $[s,t]$.

  In either case, $\xi_{u_N}$ is bounded in $L^\infty(s,t;\H)$. The bound
  $$\|P_Nf(u_N)\|_{H^{-1}}\le C\|f(u_N)\|_{L^{(2^*)'}}\le C$$
  shows that $\Dtt u_N$ is bounded in
  $L^\infty(s,t;H^{-1})$. By the Aubin--Lions--Simon theorem \cite{simon87} and the compactness
  of $L^2\hookrightarrow H^{-\de_0}$, for some $\tilde u$ and a subsequence,
  \begin{equation}\label{gal-conv}
    \begin{aligned}
       & u_N\to\tilde u\ \text{~~in  }C([s,t];H^{1-\de_0}),                                              \\
       & \Dt u_N\to\Dt\tilde u\ \text{~~in  }C([s,t];H^{-\de_0}),                                        \\
       & \xi_{u_N}(\varsigma)\wto\xi_{\tilde u}(\varsigma)\ \text{~~in  }\H\ \ \forall\varsigma\in[s,t],
    \end{aligned}
  \end{equation}
  for every $\de_0\in(0,1]$. As in Step 3 of the proof of Theorem \ref{thm-exist},
  $f(u_N)\to f(\tilde u)$ in $C([s,t];L^{\vartheta})$ for some $\vartheta>(2^*)'$, so
  $\tilde u$ is a weak solution of \eqref{eq1} on $[s,t]$. Moreover $\xi_{\tilde u}(\varsigma_0)=\xi_u(\varsigma_0)$,
  because $P_N\xi_u(\varsigma_0)\to\xi_u(\varsigma_0)$ strongly. By Corollary \ref{cor-back-unique} for $\varsigma_0=t$ and
  Theorem \ref{thm-unique} for $\varsigma_0=s$, $\tilde u=u$
  on $[s,t]$, and therefore the whole sequence converges.

  \emph{Step 3: passage to the limit.} By \eqref{gal-conv}, $u_N\to u$ in $C([s,t];L^{p+1})$.
  The functions $F(s)$ and $f(s)s$ have derivatives bounded by $C(1+|s|^p)$. H\"older's inequality together with the
  uniform $L^{p+1}$ bound therefore gives
  \[
    \int_\Om F(u_N)\to\int_\Om F(u),\qquad \la f(u_N),u_N\ra\to\la f(u),u\ra
    \qquad\text{uniformly on }[s,t].
  \]
  Next, $e^{\ga\varsigma}\(\Dt g+\frac\ga2 g\)\in L^1(s,t;H^{-1})$ and $u_N\overset{*}{\wto}u$ in
  $L^\infty(s,t;H^1_0)=\(L^1(s,t;H^{-1})\)'$. Together these give
  \begin{equation}\label{R-conv}
    \int_s^tR_N\,\td\varsigma\ \longrightarrow\ \int_s^tR\,\td\varsigma .
  \end{equation}
  At a time $\varsigma\in[s,t]$, every term of $\mathscr E_N(\varsigma)$ converges except the quadratic part
  $\frac12\|\Dt\de_N\|^2+\frac12\|\Nx\de_N\|^2$:
  \begin{itemize}
    \item $-\frac{\ga^2}8\|\de_N\|^2$ converges by \eqref{gal-conv};
    \item $e^{\ga\varsigma}\int F(u_N)$ converges by the above;
    \item $\la e^{\ga\varsigma} g(\varsigma),\de_N(\varsigma)\ra$ converges since $\de_N(\varsigma)\wto\de(\varsigma)$ in $H^1_0$.
  \end{itemize}

  The quadratic part is the squared $\H$-norm of
  $\xi_{\de_N}(\varsigma)=e^{\ga\varsigma/2}\(u_N,\Dt u_N+\frac\ga2u_N\)(\varsigma)$, which converges weakly by
  \eqref{gal-conv}. Hence
  \begin{equation}\label{E-lsc}
    \mathscr E(\varsigma)\le\liminf_{N\to\infty}\mathscr E_N(\varsigma)\quad\forall\varsigma\in[s,t],\qquad
    \mathscr E_N(\varsigma_0)\to\mathscr E(\varsigma_0),
  \end{equation}
  the second statement because the data at $\varsigma_0$ converge strongly in $\H$.

  \emph{Step 4: two inequalities.} Taking $\varsigma_0=s$ and using \eqref{gal-exp}, \eqref{R-conv},
  \eqref{E-lsc},
  \[
    \mathscr E(t)\le\liminf_N\mathscr E_N(t)
    =\lim_N\Big(\mathscr E_N(s)-\int_s^tR_N\Big)=\mathscr E(s)-\int_s^tR .
  \]
  Taking $\varsigma_0=t$ in the same way,
  \[
    \mathscr E(s)\le\liminf_N\mathscr E_N(s)
    =\lim_N\Big(\mathscr E_N(t)+\int_s^tR_N\Big)=\mathscr E(t)+\int_s^tR .
  \]
  Together these give \eqref{expid}.
\end{proof}

\begin{cor}[Strong continuity and the energy identity]\label{cor-energy-eq}
  Under the assumptions of Lemma \ref{lem-expenergy}, let $u$ be a weak solution of \eqref{eq1} on $[s_0,t_0]$.
  Then $\xi_u\in C([s_0,t_0];\H)$. Moreover, the energy equality holds,
  \begin{equation}\label{energy-eq}
    E\(u(t)\)+\ga\int_s^t\|\Dt u\|^2\,\td\varsigma=E\(u(s)\)-\int_s^t\la\Dt g,u\ra\,\td\varsigma,
    \qquad s_0\le s\le t\le t_0 .
  \end{equation}
\end{cor}

\begin{proof}
  \emph{Strong continuity.} By \eqref{expid}, $\mathscr E$ is continuous on $[s_0,t_0]$.
  All its terms other than the quadratic part are continuous on $[s_0,t_0]$:
  \begin{itemize}
    \item $\|\de\|^2$, since $u\in C(L^2)$;
    \item $e^{\ga\varsigma}\int F(u)$, since $u\in C(L^\eta)$ with $\eta\in[p+1,2^*)$;
    \item $\la g(\varsigma),\de\ra$, which pairs $g\in C(H^{-1})$ (strongly continuous) with the bounded, weakly continuous $\de\in C_w(H^1_0)$.
  \end{itemize}

  Hence $\varsigma\mapsto\|\xi_\de(\varsigma)\|_{\H}$ is continuous. Since $\xi_\de$ is weakly continuous
  in the Hilbert space $\H$, it is strongly continuous, and so is
  $\xi_u=e^{-\ga\varsigma/2}\(\de,\Dt\de-\frac\ga2\de\)$.

  \emph{Energy identity.} Fix $s<t$ and let $u_N$ be the Galerkin solution with data
  $P_N\xi_u(s)$ at $s$. For every $\varsigma\in[s,t]$, \eqref{R-conv} and \eqref{expid} on $[s,\varsigma]$
  give the full limit
  $$\mathscr E_N(\varsigma)\to\mathscr E(s)-\int_s^{\varsigma}R=\mathscr E(\varsigma).$$
  By Step 3 above, this is equivalent to $\|\xi_{\de_N}(\varsigma)\|_{\H}\to\|\xi_\de(\varsigma)\|_{\H}$. Together with weak
  convergence, $\xi_{u_N}(\varsigma)\to\xi_u(\varsigma)$ strongly in $\H$ for every $\varsigma\in[s,t]$. Since
  $\xi_{u_N}$ is bounded, dominated convergence gives $\Dt u_N\to\Dt u$ in $L^2(s,t;L^2(\Om))$.
  Multiplying the Galerkin equation by $\Dt u_N$ and integrating, the Galerkin solutions satisfy
  \[
    E\(u_N(t)\)+\ga\int_s^t\|\Dt u_N\|^2=E\(u_N(s)\)-\int_s^t\la\Dt g,u_N\ra .
  \]
  Every term now converges: the endpoint energies by strong convergence in $\H$ and
  $p+1<2^*$, the dissipation integral by the $L^2_tL^2_x$ convergence, and the force term by
  weak-$*$ convergence. This gives \eqref{energy-eq}.
\end{proof}

\section{Pullback attractors}\label{sec-attr}

By Theorems \ref{thm-unique}, \ref{thm-exist} and \ref{thm-lip}, the problem
\begin{equation}\label{eq-g}
  \Dtt u+\ga\Dt u-\Dx u+f(u)=g(t),\qquad u|_{\partial\Om}=0,\qquad
  \xi_u|_{t=\tau}=\xi_\tau\in\H
\end{equation}
generates a \emph{process} $U(t,\tau):\H\to\H$, $U(t,\tau)\xi_\tau:=\xi_u(t)$, $t\ge\tau$, in
the sense of Chapter 1 of \cite{clr13}: $U(t,s)U(s,\tau)=U(t,\tau)$ for $\tau\le s\le t$,
$U(\tau,\tau)=\mathrm{Id}$, and $(\tau,t,\xi_\tau)\mapsto U(t,\tau)\xi_\tau$ is continuous,
the last property following from $u\in C([\tau,T];H^1_0)\cap C^1([\tau,T];L^2)$ in Definition
\ref{def-weak} together with the Lipschitz estimate of Theorem \ref{thm-lip}. Throughout this
section the force $g$ is the single, fixed function of hypothesis \eqref{G1}; no translate or
approximant of $g$ is ever introduced.

\subsection{Dissipativity}\label{ssec-pb-dissip}

By \eqref{extra-reg},
\begin{equation}\label{dissip-single}
  \|U(t,\tau)\xi_\tau\|^2_{\H}\le Q\bigl(\|\xi_\tau\|_{\H}+\|g\|_{H^1_b}\bigr)e^{-\al(t-\tau)}+Q_0,
  \qquad Q_0:=Q\bigl(\|g\|_{H^1_b}\bigr),\qquad t\ge\tau.
\end{equation}
Fix $t\in\R$ and let $D\subset\H$ be bounded, say $D\subset\{\|\xi\|_{\H}\le R\}$. For $\tau\le t$,
\eqref{dissip-single} gives
\[
  \sup_{\xi_\tau\in D}\|U(t,\tau)\xi_\tau\|^2_{\H}\le Q(R)e^{-\al(t-\tau)}+Q_0,
\]
whose right-hand side converges to $Q_0$ as $\tau\to-\infty$. The ball
\begin{equation}\label{B0}
  \mathcal B_0:=\big\{\xi\in\H:\ \|\xi\|_{\H}^2\le 2Q_0+1\big\}
\end{equation}
therefore pullback attracts $D$ at time $t$: for $\tau$ sufficiently negative one has
$U(t,\tau)D\subset\mathcal B_0$. Since $\mathcal B_0$ is independent of $t$, the same holds at every time
$\varsigma\in\R$, so $U(\cdot,\cdot)$ is strongly pullback bounded dissipative in the sense of
Definition \ref{def-pb-dissip}, the single ball $\mathcal B_0$ serving at every time.

\subsection{Pullback asymptotic compactness}\label{ssec-pb-compact}

This is the only nontrivial hypothesis of Theorem \ref{thm-CLR-existence}, and we verify it by the energy
method, reusing the exponentially weighted identity of Lemma \ref{lem-expenergy} rather than reproving it. Working
with a single, fixed force makes this reuse literal: Lemma \ref{lem-expenergy} was proved for Galerkin
approximations of \eqref{eq1} with the force $g$ fixed throughout, whereas the sequence constructed below consists
of \emph{genuine} weak solutions with that same $g$. The identity \eqref{expid} therefore applies to each of them
directly, with equality, and none of the approximation machinery of Section \ref{ssec-expenergy} needs to be
repeated.

\begin{thm}\label{thm-pb-compact}
  Assume \eqref{F1}, \eqref{F2} and \eqref{G1}. Then $U(\cdot,\cdot)$ is pullback
  asymptotically compact.
\end{thm}

\begin{proof}
  Fix $t\in\R$, let $s_k\to-\infty$ and let $\{x_k\}\subset\H$ be bounded, say
  $\|x_k\|_{\H}\le R$. Let $u_k$ be the weak solution of \eqref{eq-g} with $\xi_{u_k}(s_k)=x_k$,
  given globally by Theorem \ref{thm-exist}, and set $\xi_k(s):=\xi_{u_k}(s)$ for
  $s\ge s_k$.

  \emph{Step 1: diagonal extraction and limit trajectory.} By \eqref{dissip-single},
  $$\sup_{s\ge s_k}\|\xi_k(s)\|_{\H}^2\le Q(R+\|g\|_{H^1_b})+Q_0$$
  for every $k$, uniformly.

  For each integer $j\ge1$, consider the window $I_j:=[t-j,t]$. For $k$ large enough that $s_k<t-j$,
  $u_k$ solves \eqref{eq-g} on $I_j$ with $\sup_k\|\xi_k\|_{L^\infty(I_j;\H)}<\infty$. As in Step 3
  of Theorem \ref{thm-exist}, $\Dtt u_k$ is bounded in $L^\infty(I_j;H^{-1})$, and the
  Aubin--Lions--Simon theorem \cite{simon87} gives a subsequence (call it $\{u_k^{(j)}\}$) with
  \begin{equation}\label{pb-conv-j}
    \xi_k^{(j)}\to \xi^{(j)}\ \text{in}\ C(I_j;H^{1-\de}(\Om)\times H^{-\de}(\Om))\ \ \forall\de\in(0,1],\qquad
    \xi_k^{(j)}\overset{*}{\wto}\xi_{u^{(j)}}\ \text{in}\ L^\infty(I_j;\H),
  \end{equation}
  and $f(u_k^{(j)})\to f(u^{(j)})$ in $C(I_j;L^{\vartheta}(\Om))$ for some $\vartheta>(2^*)'$,
  so $u^{(j)}$ is a weak solution of \eqref{eq-g} on $I_j$.

  By a diagonal argument we now extract a single subsequence, still denoted $\{u_k\}$, which is a subsequence of
  every $\{u_k^{(j)}\}$. Then on each window $I_j$ the convergence above holds with the limit $u^{(j)}$.

  We claim that these limits agree on the overlaps. Indeed, fix $j$ and let $\zeta\in[t-j-1,t-j]$ be a point common
  to $I_j$ and $I_{j+1}$. On $I_j$ the extracted subsequence satisfies $\xi_k(\zeta)\wto\xi^{(j)}(\zeta)$ in $\H$,
  and on $I_{j+1}$ it satisfies $\xi_k(\zeta)\wto\xi^{(j+1)}(\zeta)$; since the weak limit is unique, the two coincide. Hence the pieces $u^{(j)}$ assemble into a single function
  $u:(-\infty,t]\to H^1_0(\Om)$ with $u|_{I_j}=u^{(j)}$, which is a weak solution on every $I_j$, and therefore
  on $(-\infty,t]$. By weak lower semicontinuity of the $\H$-norm,
  \[
    \sup_{s\le t}\|\xi_u(s)\|^2_{\H}\le\liminf_k\sup_{s\le t}\|\xi_k(s)\|^2_{\H}
    \le Q(R+\|g\|_{H^1_b})+Q_0<\infty,
  \]
  the last inequality being \eqref{dissip-single}.

  In particular, on every window $I_j$,
  \begin{equation}\label{pb-conv}
    \xi_k\to \xi_u\ \text{in}\ C(I_j;H^{1-\de}(\Om)\times H^{-\de}(\Om))\ \ \forall\de\in(0,1],\qquad
    \xi_k\overset{*}{\wto}\xi_u\ \text{in}\ L^\infty(I_j;\H),
  \end{equation}
  so that for every $s\le t$ we have $u_k(s)\to u(s)$ in $H^{1-\de}(\Om)$ for all $\de>0$, and
  $\xi_k(s)\wto\xi_u(s)$ in $\H$.

  \emph{Step 2: passage to the limit in the weighted energy identity.}
  Apply the identity \eqref{expid} of Lemma \ref{lem-expenergy} to each $u_k$ on $[s_k,t]$: since
  $u_k$ is a genuine weak solution (not merely a Galerkin approximation), the identity holds
  with equality,
  \begin{equation}\label{pb-exact}
    \mathscr E_k(t)=\mathscr E_k(s_k)-\int_{s_k}^tR_k(s)\,\td s,
  \end{equation}
  where
  \begin{align*}
    \mathscr E_k(t) & =e^{\ga t}\(\frac12\|\Dt u_k(t)\|^2+\frac12\|\Nx u_k(t)\|^2
    +\frac{\ga}{2}\la \Dt u_k(t),u_k(t)\ra+\!\int_\Om F\(u_k(t)\)\td x
    -\la g(t),u_k(t)\ra\),                                                                \\
    R_k(s)          & =\ga e^{\ga s}\!\int_\Om\Big[\frac12f(u_k)u_k-F(u_k)\Big](s)\,\td x
    +e^{\ga s}\Big\la\Dt g(s)+\frac\ga2g(s),u_k(s)\Big\ra.
  \end{align*}

  We let $k\to\infty$
  in \eqref{pb-exact}, i.e. $s_k\to-\infty$.

  The boundary term $\mathscr E_k(s_k)$ is a sum of terms each bounded by a
  fixed monotone function of $\|\xi_k(s_k)\|_{\H}=\|x_k\|_{\H}\le R$ and of
  $\|g(s_k)\|_{H^{-1}}\le\sup_{s\in\R}\|g(s)\|_{H^{-1}}<\infty$ (by \eqref{gsup}), multiplied by
  $e^{\ga s_k}$; hence
  \begin{equation}\label{Ek-limit}
    |\mathscr E_k(s_k)|\le C(R)e^{\ga s_k}\to0.
  \end{equation}

  For the integral term, we claim $\mathbf 1_{[s_k,t]}(s)R_k(s)\to R(s)$ for a.e.\ $s<t$. Indeed for any $s\le t$, Step 1 gives
  $$u_k(s)\to u(s) \quad\text{in}\quad L^{p+1}(\Om) \text{\quad and \quad}  u_k(s) \wto u(s)\quad\text{weakly in }H^1_0(\Om).$$
  By \eqref{fbounds},
  \[
    \left|\frac{1}{2}f(u)u - F(u)\right| \le C(1 + |u|^{p+1}).
  \]
  Therefore,
  \[
    \int_\Om\Big[\frac{1}{2}f(u_k(s))u_k(s) - F(u_k(s))\Big]\td x
    \to \int_\Om\Big[\frac{1}{2}f(u(s))u(s) - F(u(s))\Big]\td x
  \]
  for almost every $s \in (-\infty, t]$. Since $u_k(s)\wto u(s)$ weakly in $H^1_0(\Om)$ and $\Dt g(s) + \frac{\ga}{2}g(s)\in H^{-1}(\Om)$ is fixed,
  the duality pairing converges:
  $$\Big\langle\Dt g(s) + \frac{\ga }{2}g(s), u_k(s)\Big\rangle\to \Big\langle\Dt g(s) + \frac{\ga}{2}g(s), u(s)\Big\rangle.$$
  Combining these, $\mathbf 1_{[s_k,t]}(s)R_k(s)\to R(s)$ for a.e.\ $s\in(-\infty,t)$.

  To apply dominated convergence to $\int_{s_k}^tR_k(s)\,\td s$, note that
  \begin{align*}
    \ga e^{\ga s} \int_\Om \left|\frac{1}{2}f(u_k(s))u_k(s) - F(u_k(s))\right|\,\td x
     & \le C\ga e^{\ga s}\big(|\Om| + \|u_k(s)\|_{L^{p+1}(\Om)}^{p+1}\big) \\
     & \le C(R)\ga e^{\ga s},
  \end{align*}
  and the uniform bound $\|\xi_k\|_{L^\infty(s_k,t;\H)}\le C(R)$,
  \[
    \left|e^{\ga s}\Big\la\Dt g(s)+\frac\ga2g(s),u_k(s)\Big\ra\right|\le C(R)e^{\ga s}\big(1+\|\Dt g(s)\|_{H^{-1}}+\|g(s)\|_{H^{-1}}\big).
  \]
  Since $g\in H^1_b(\R;H^{-1})$, the right-hand side belongs to $L^1(-\infty,t)$:
  \begin{equation}
    \begin{aligned}
      \int_{-\infty}^te^{\ga s}\left(\|\Dt g(s)\|_{H^{-1}}+\|g(s)\|_{H^{-1}}\right)\,\td s
       & =\sum_{j=0}^{\infty}\int_{t-j-1}^{t-j} e^{\ga s}\left(\|\Dt g(s)\|_{H^{-1}}+\|g(s)\|_{H^{-1}}\right)\,\td s \\
       & \le \sum_{j=0}^{\infty}e^{\ga (t-j)}\|g\|_{L^1(t-j-1,t-j;H^{-1})}                                           \\
       & \le \|g\|_{H^1_b}\sum_{j=0}^{\infty}e^{\ga (t-j)}                                                           \\
       & \le Ce^{\ga t}\|g\|_{H^1_b}<\infty.
    \end{aligned}
  \end{equation}
  Therefore dominated convergence gives
  \begin{equation}\label{Rk-limit}
    \int_{s_k}^tR_k(s)\,\td s=\int_{-\infty}^{t}\mathbf 1_{[s_k,t]}(s)R_k(s)\,\td s\to\int_{-\infty}^tR(s)\,\td s.
  \end{equation}

  Finally, since $u_k(t)\to u(t)$ strongly in $L^{p+1}(\Om)$ and
  weakly in $H^1_0(\Om)$, while $\Dt u_k(t)\rightharpoonup\Dt u(t)$ weakly in $L^2(\Om)$.
  It follows that in the term $\mathscr E_k(t)$,
  \[
    \int_\Om F(u_k(t))\,\td x\to\int_\Om F(u(t))\,\td x,
    \qquad
    \la g(t),u_k(t)\ra\to\la g(t),u(t)\ra.
  \]
  The cross term also converges, since $u_k(t)\to u(t)$ strongly in $L^2(\Om)$ and
  $\Dt u_k(t)\rightharpoonup\Dt u(t)$ weakly in $L^2(\Om)$:
  \[
    \la \Dt u_k(t),u_k(t)\ra\longrightarrow
    \la \Dt u(t),u(t)\ra.
  \]
  Using these convergences and lower semicontinuity of $\H$-norm, we have
  \begin{equation}\label{Ek-lsc}
    \mathscr E(t) \le \liminf_{k\to\infty} \mathscr E_k(t).
  \end{equation}

  Substituting \eqref{Ek-limit}, \eqref{Rk-limit} and \eqref{Ek-lsc} into \eqref{pb-exact} yields
  \begin{equation}\label{Ek-liminf}
    \mathscr E(t)\le \liminf_{k\to \infty}\mathscr E_k(t)=-\int_{-\infty}^t R(s)\,\td s,
  \end{equation}

  \emph{Step 3: strong convergence at the terminal time.}
  On the other hand, the same identity \eqref{expid} applied to $u$ gives
  \begin{equation}
    \mathscr E(t)=\mathscr E(s)-\int_{s}^t R(\varsigma)\,\td\varsigma .
  \end{equation}
  The boundary term $\mathscr E(s)$ tends to $0$ as $s\to-\infty$, by the same estimate that controlled
  $\mathscr E_k(s_k)$ and using
  $\sup_{\varsigma\le t}\|\xi_u(\varsigma)\|^2_{\H}\le Q(R+\|g\|_{H^1_b})+Q_0<\infty$. Hence
  \begin{equation}\label{integral-R}
    \mathscr E(t)=-\int_{-\infty}^tR(\varsigma)\,\td\varsigma,
  \end{equation}
  Comparing with \eqref{Ek-liminf} yields
  \[
    \mathscr E(t)=\liminf_{k\to\infty}\mathscr E_k(t).
  \]
  Since every other term of $\mathscr E_k(t)$ converges, this is equivalent to
  $\|\xi_k(t)\|_{\H}\to\|\xi_u(t)\|_{\H}$. Together with the weak convergence
  $\xi_k(t)\wto\xi_u(t)$ in the Hilbert space $\H$, it gives the strong convergence
  $\xi_k(t)\to\xi_u(t)$ in $\H$, which is the convergent subsequence we needed.
\end{proof}

\subsection{Existence of the pullback attractor}\label{ssec-pb-exist}

\begin{thm}\label{thm-pullback}
  Assume \eqref{F1}, \eqref{F2} and \eqref{G1}. Then the process $U(\cdot,\cdot)$ generated by
  the weak solutions of \eqref{eq1} possesses a pullback attractor $\A(\cdot)\subset\H$,
  bounded in the past, with
  \begin{equation}\label{pullback-rep}
    \A(t)=\big\{u(t):u\ \text{a complete, backwards-bounded weak solution of \eqref{eq1}}\big\},
    \qquad t\in\R .
  \end{equation}
\end{thm}

\begin{proof}
  By Section \ref{ssec-pb-dissip}, $U(\cdot,\cdot)$ is strongly pullback bounded dissipative, with
  the single ball $\mathcal B_0$ of \eqref{B0} pullback attracting bounded sets at every time.
  By Theorem \ref{thm-pb-compact}, $U(\cdot,\cdot)$ is pullback asymptotically compact. Theorem
  \ref{thm-CLR-existence} then yields a pullback attractor $\A(\cdot)$, bounded in the past.
  Theorem \ref{thm-CLR-kernel} identifies the attractor as in \eqref{pullback-rep}: a global solution
  $\varphi$ with $\sup_{s\le t}\|\varphi(s)\|_{\H}<\infty$ is precisely a weak solution $u$ of \eqref{eq1} defined
  on all of $(-\infty,t]$ and bounded backwards from $t$, that is, a complete solution bounded in the past.
\end{proof}

\begin{remark}
  Theorem \ref{thm-pullback} neither requires nor produces an attraction property uniform in the terminal time
  $t$: for each fixed $t$, the set $\A(t)$ attracts bounded sets as the initial time is pulled back to $-\infty$,
  but both the rate of attraction and the set $\A(t)$ itself may depend on $t$. Note, however, that by the
  time-independence of $\mathcal B_0$ established in Section \ref{ssec-pb-dissip}, the union
  $\bigcup_{t\in\R}\A(t)$ is in fact bounded in $\H$. This is both the price and the simplification of working with
  the single process generated by $g$ rather than with the family indexed by its hull.
\end{remark}

\appendix
\section{The optimal exponents within the duality scheme}\label{app-opt}

This appendix determines the largest exponent permitted by the present duality scheme. It should not be read as
a claim that the underlying bounded-domain Strichartz range is globally optimal; indeed, the admissible range of
\cite{blair09} is not known to be optimal in general. Throughout, $n\ge3$, $I$ is an interval with $|I|\le1$, and
admissibility is understood in the sense of Definition \ref{def-adm}. We optimise over \emph{all} regularity indices
$\sg\in[0,1]$ simultaneously; it will turn out that the largest growth exponent does not depend on $\sg$ at all ---
a fact that is invisible if one fixes $\sg=1$ from the start, and that explains why the whole-space theory of
\cite{glassey82} is forced to work at $\sg=\frac12$, whereas we may use the endpoint choice $\sg=1$.

\begin{Def}\label{def-scheme}
  Let $1<p<p_c$. A \emph{duality scheme for the exponent $p$} on $\Om$ consists of an index
  $\sg\in[0,1]$ and two triples $(q,r,\sg)$, $(\tilde q,\tilde r,1-\sg)$, both
  $\Om$-admissible, such that
  \begin{equation}\label{scheme-cond}
    \frac1{\tilde r'}=\frac1\mu+\frac1r,\qquad \mu=\frac{2^*}{p-1},
    \qquad\text{and}\qquad
    \theta:=\frac1{\tilde q'}-\frac1q>0.
  \end{equation}
  A duality scheme on $\R^n$ is defined in the same way with $\R^n$-admissibility.
\end{Def}

For any $u, v\in H^1_0(\Omega)$,
\begin{equation}
  f(u)-f(v)=V\,(u-v),\qquad V:=\int_0^1f'\(su+(1-s)v\)\,\td s.
\end{equation}
The first condition in \eqref{scheme-cond} is exactly what makes the H\"older estimate
\begin{equation}\label{holder-scheme}
  \big\|V(u-v)\big\|_{L^{\tilde q'}(I;L^{\tilde r'}(\Om))}
  \le|I|^{\theta}\,\|V\|_{L^\infty(I;L^{\mu}(\Om))}\,\|u-v\|_{L^{q}(I;L^{r}(\Om))}
\end{equation}
available. The second produces the positive power of $|I|$ without which the nonlinear
term cannot be absorbed on short intervals.

We use the variables
\begin{equation}\label{vars}
  a=\frac1q,\quad b=\frac1r,\quad c=\frac1{\tilde q},\quad d=\frac1{\tilde r},
  \qquad S:=b+d,
\end{equation}
and we write the two boundary conditions of Definition \ref{def-adm} in the unified form
\begin{equation}\label{unified}
  \kappa a+\ell b\le L,\qquad \kappa c+\ell d\le L,
\end{equation}
with the parameters
\begin{equation}\label{kappa-ell}
  \begin{array}{c|ccc}
                & \kappa & \ell & L             \\\hline
    \Om,\ n=3   & 3      & 2    & 1             \\[1mm]
    \Om,\ n\ge4 & 1      & 1    & \frac12       \\[1mm]
    \R^n        & 2      & n-1  & \frac{n-1}{2}
  \end{array}
\end{equation}
Here the first row is the case $n=3$ of \eqref{adm-om}; for $n=4$ that
branch coincides with the second row, both giving $\frac1q+\frac1r\le\frac12$, so no
generality is lost by listing only $n=3$ separately.
In all three cases
\begin{equation}\label{L-half}
  L=\frac{\ell}{2},
\end{equation}
a coincidence which is responsible for the clean form of everything below. Finally set
\begin{equation}\label{M-def}
  M:=\kappa n-\ell>0 .
\end{equation}

\begin{lemma}\label{lem-alg}
  Let $(q,r,\sg)$ and $(\tilde q,\tilde r,1-\sg)$ satisfy the two scaling identities of
  Definition \ref{def-adm}. Then
  \begin{equation}\label{alg1}
    (a+c)+nS=n-1 .
  \end{equation}
  If in addition the first condition in \eqref{scheme-cond} holds, then
  \begin{equation}\label{alg2}
    p=1+2^*\(1-S\),\qquad\text{equivalently}\qquad S=1-\frac{p-1}{2^*},
  \end{equation}
  and
  \begin{equation}\label{alg3}
    \theta=1-(a+c)=nS-n+2=2-\frac{(n-2)(p-1)}{2}.
  \end{equation}
\end{lemma}

\begin{proof}
  The scaling identities read $a+nb=\frac n2-\sg$ and $c+nd=\frac n2-(1-\sg)$; adding them
  eliminates $\sg$ and gives \eqref{alg1}. The first condition in \eqref{scheme-cond} is
  $1-d=\frac1\mu+b$, i.e. $\frac1\mu=1-S$, and since $\frac1\mu=\frac{p-1}{2^*}$ this is
  \eqref{alg2}. Finally $\theta=(1-c)-a=1-(a+c)$, which by \eqref{alg1} equals $nS-n+2$;
  substituting \eqref{alg2} and $2^*=\frac{2n}{n-2}$ gives
  $\theta=n\big(1-\frac{(n-2)(p-1)}{2n}\big)-n+2$, as claimed.
\end{proof}

Note that $\sg$ has already disappeared from \eqref{alg1}--\eqref{alg3}.

\begin{thm}\label{thm-opt}
  Let $n\ge3$ and let $\kappa,\ell,L,M$ be as in \eqref{kappa-ell}--\eqref{M-def}.
  \begin{enumerate}
    \item\label{opt-S} Every duality scheme satisfies
    \begin{equation}\label{Smin}
      S\ \ge\ S_{\min}:=\frac{\kappa(n-1)-\ell}{M},
    \end{equation}
    and consequently
    \begin{equation}\label{pmax}
      p\ \le\ p_{\max}:=1+2^*\cdot\frac{\kappa}{M}
      =1+\frac{2n}{n-2}\cdot\frac{\kappa}{\kappa n-\ell}.
    \end{equation}
    Both bounds are independent of $\sg$.
    \item\label{opt-family} Equality in \eqref{Smin} holds precisely when both inequalities in
    \eqref{unified} are equalities, which happens exactly for
    \begin{equation}\label{optfamily}
      a=\frac{\sg\ell}{M},\quad b=\frac12-\frac{\kappa\sg}{M},\qquad
      c=\frac{(1-\sg)\ell}{M},\quad d=\frac12-\frac{\kappa(1-\sg)}{M},
    \end{equation}
    a one-parameter family indexed by $\sg$.
    \item\label{opt-window} The member \eqref{optfamily} is admissible if and only if
    \begin{equation}\label{window}
      1-\sg^*<\sg<\sg^*,\qquad \sg^*:=\frac{M}{2\max\{\kappa,\ell\}},
    \end{equation}
    where the strict inequalities come from $q,\tilde q>2$.
    \item\label{opt-om} For $\Om$, one has $\sg^*>1$ in every dimension $n\ge3$, so $\sg=1$ is
    admissible; \eqref{pmax} gives
    \[
      p_{\max}=p_\Om(n)=
      \begin{cases}
        \dfrac{25}{7},            & n=3,   \\[3mm]
        1+\dfrac{2n}{(n-2)(n-1)}, & n\ge4,
      \end{cases}
    \]
    and at
    $\sg=1$ the family \eqref{optfamily} reduces to $c=0$, $d=\frac12$, i.e.
    $(\tilde q,\tilde r)=(\infty,2)$, and to the direct pair
    \begin{equation}
      (q,r)=
      \begin{cases}
        \Big(\dfrac{7}{2},\ 14\Big),         & n=3,   \\[3mm]
        \Big(n-1,\ \dfrac{2(n-1)}{n-3}\Big), & n\ge4.
      \end{cases}
    \end{equation}
    \item\label{opt-free} For $\R^n$,  one has $\sg^*\le1$, with equality only when $n=3$; hence
    $\sg=1$ is never admissible. Here \eqref{pmax} gives
    \[
      p_{\max}=1+\frac{4n}{(n-2)(n+1)}=\frac{n^2+3n-2}{n^2-n-2}=p_1(n),
    \]
    the exponent appearing in Theorem 1 of \cite{glassey82}, and the symmetric member $\sg=\frac12$ of
    \eqref{optfamily} is the conformal self-dual triple
    $q=\tilde q=r=\tilde r=\frac{2(n+1)}{n-1}$.
  \end{enumerate}
\end{thm}

\begin{proof}
  \emph{(i).} Adding the two inequalities \eqref{unified} gives
  $\kappa(a+c)+\ell S\le 2L=\ell$. Substituting
  $a+c=n-1-nS$ from \eqref{alg1},
  \[
    \kappa(n-1)-\kappa nS+\ell S\le\ell
    \quad\Longleftrightarrow\quad
    \kappa(n-1)-\ell\le(\kappa n-\ell)S,
  \]
  which is \eqref{Smin}. Then \eqref{alg2} and
  \[
    1-S_{\min}=\frac{M-\kappa(n-1)+\ell}{M}=\frac{\kappa n-\ell-\kappa n+\kappa+\ell}{M}
    =\frac{\kappa}{M}
  \]
  give \eqref{pmax}. Neither computation involves $\sg$.

  \emph{(ii).} Equality in the sum of \eqref{unified} forces equality in each term, so that
  \begin{equation}\label{eq-equality}
    \kappa a+\ell b=\frac{\ell}{2},\qquad \kappa c+\ell d=\frac{\ell}{2}.
  \end{equation}
  We solve the first of these together with the scaling identity $a+nb=\frac n2-\sg$. From
  \eqref{eq-equality},
  \begin{equation}\label{eq-aofb}
    a=\frac{\ell}{\kappa}\Big(\frac12-b\Big),
  \end{equation}
  and substituting this into the scaling identity gives
  $$\frac{\ell}{\kappa}\Big(\frac12-b\Big)+nb=\frac n2-\sg,$$
  which yields
  $$b=\frac12-\frac{\kappa\sg}{M}.$$
  And \eqref{eq-aofb} returns
  $$a=\frac{\ell}{\kappa}\cdot\frac{\kappa\sg}{M}=\frac{\sg\ell}{M}.$$
  The identities for
  $c,d$ follow from the second equation in \eqref{eq-equality} in the same way, with $\sg$
  replaced by $1-\sg$.

  \emph{(iii).} The constraints defining admissibility are $0<a<\frac12$ and $0<b\le\frac12$,
  and likewise for $c,d$. By \eqref{optfamily} the binding ones are
  \[
    b>0\iff\sg<\frac{M}{2\kappa},\qquad
    a<\frac12\iff\sg<\frac{M}{2\ell},
  \]
  together with the two conditions obtained by replacing $\sg$ with $1-\sg$. Combining the four
  gives \eqref{window}.

  \emph{(iv) and (v).} By \eqref{pmax} and \eqref{window} everything is determined by the three
  numbers $\kappa,\ell,M$, which we now tabulate:
  \begin{equation}\label{geomtable}
    \begin{array}{c|ccc|cc}
                  & \kappa                          & \ell & M   & \sg^*=\dfrac{M}{2\max\{\kappa,\ell\}}
                  & p_{\max}=1+2^*\dfrac{\kappa}{M}                                                                      \\\hline
                  &                                 &      &     &                                       &               \\[-2mm]
      \Om,\ n=3   & 3                               & 2    & 7   & \dfrac76                              & \dfrac{25}{7} \\[3mm]
      \Om,\ n\ge4 & 1                               & 1    & n-1 & \dfrac{n-1}{2}
                  & 1+\dfrac{2n}{(n-2)(n-1)}                                                                             \\[3mm]
      \R^n        & 2                               & n-1  & n+1 & \dfrac{n+1}{2(n-1)}                   & p_1(n)
    \end{array}
  \end{equation}

  This proves the assertions about $\sg$. In the two cases on $\Om$ we have $\sg^*>1$, so
  $\sg=1$ satisfies \eqref{window}, the lower constraint being vacuous because $1-\sg^*<0$. On
  $\R^n$ we have $\sg^*\le1$, so $\sg=1$ violates the strict upper inequality in
  \eqref{window} and is never admissible, whereas $\sg=\frac12$ is admissible in every
  dimension since $\sg^*>\frac12$.

  The value of $p_{\max}$ in the last column is \eqref{pmax} evaluated on the table. For
  $\R^n$ this reads
  \begin{equation}\label{p1-identify}
    p_{\max}=1+\frac{2n}{n-2}\cdot\frac{2}{n+1}
    =\frac{(n-2)(n+1)+4n}{(n-2)(n+1)}
    =\frac{n^2+3n-2}{n^2-n-2}=p_1(n),
  \end{equation}
  which identifies the threshold with the exponent appearing in Theorem 1 of \cite{glassey82}.

  It remains to compute the corresponding exponent pairs at the endpoint of this scheme. On $\Om$ we take $\sg=1$ in
  \eqref{optfamily}, which gives
  \begin{equation}\label{sigma1-red}
    c=0,\quad d=\frac12,\qquad a=\frac{\ell}{M},\quad b=\frac12-\frac{\kappa}{M},
  \end{equation}
  so that $(\tilde q,\tilde r)=(\infty,2)$ and the source norm is $L^1(I;L^2(\Om))$.
  Numerically, $a=\frac27$ and $b=\frac1{14}$ for $n=3$, whereas $a=\frac1{n-1}$ and $b=\frac{n-3}{2(n-1)}$ for $n\ge4$; inverting yields
  the pairs.

  On $\R^n$ we take instead the symmetric value $\sg=\frac12$, for which \eqref{optfamily} gives
  \[
    a=c=\frac{n-1}{2(n+1)},\qquad b=d=\frac12-\frac{1}{n+1}=\frac{n-1}{2(n+1)} .
  \]
  All four exponents coincide, so $q=\tilde q=r=\tilde r=\frac{2(n+1)}{n-1}$ and the scheme is
  self-dual.
\end{proof}

\begin{remark}[The role of $\sg$]\label{rem-sigma}
  Theorem \ref{thm-opt}\eqref{opt-family} shows that the optimum is not a single point but a curve parametrised by
  $\sg$, along which $p_{\max}$, $S_{\min}$ and $\theta$ are all constant. Thus $\sg$ is not an optimisation
  variable at all: it merely selects a point on the endpoint curve of this scheme, and its only genuine effect is
  through the admissibility window \eqref{window}. This is where the two geometries part company. On a bounded domain
  the window contains $\sg=1$, and that choice is the most convenient one, since by \eqref{optfamily} it collapses
  the dual triple to $(\infty,2)$ and the source norm to $L^1(I;L^2(\Om))$, the simplest possible. In the whole
  space $\sg=1$ is excluded, and \cite{glassey82} is obliged to work at the symmetric point
  $\sg=\frac12$, where the scheme is self-dual and the source norm is the conformal
  $L^{4/3}_{t,x}$ in dimension three.
\end{remark}



\begin{thebibliography}{00}

  \bibitem{babin92} A.V. Babin, M.I. Vishik, \emph{Attractors of evolution equations},
  North-Holland, Amsterdam, 1992.

  \bibitem{ball04} J.M. Ball, Global attractors for damped semilinear wave equations,
  \emph{Discrete Contin. Dyn. Syst.}, \textbf{10} (2004) 31--52.

  \bibitem{blair09} M.D. Blair, H.F. Smith, C.D. Sogge, Strichartz estimates for the wave
  equation on manifolds with boundary, \emph{Ann. Inst. H. Poincar\'e C Anal. Non Lin\'eaire},
  \textbf{26} (2009) 1817--1829.

  \bibitem{bulut} A. Bulut, M. Czubak, D. Li, N. Pavlovi\'c, X. Zhang, Stability and Unconditional
  Uniqueness of Solutions for Energy Critical Wave Equations in High Dimensions,
  \emph{Communications in Partial Differential Equations}, \textbf{38} (2013) 575--607.

  \bibitem{burq} N. Burq, G. Lebeau, F. Planchon, Global existence for energy critical waves in $3$-d domains,
  \emph{Journal of the American Mathematical Society}, \textbf{21} (2008) 831--845.

  \bibitem{cara04} T. Caraballo, P. E. Kloeden, J. Real, Pullback and forward attractors for a
  damped wave equation with delays, \emph{Stoch. Dyn.}, \textbf{4} (2004), 405--423.

  \bibitem{clr13} A.N. Carvalho, J.A. Langa, J.C. Robinson, \emph{Attractors for
    infinite-dimensional non-autonomous dynamical systems}, Applied Mathematical Sciences,
  vol.~182, Springer, New York, 2012.

  \bibitem{chueshov02} V.V. Chepyzhov, M.I. Vishik, \emph{Attractors for equations of
    mathematical physics}, AMS Colloquium Publications, vol.~49, American Mathematical Society,
  Providence, RI, 2002.



  \bibitem{glassey82} R.T. Glassey, M. Tsutsumi, On uniqueness of weak solutions to
  semilinear wave equations, \emph{Comm. Partial Differential Equations}, \textbf{7} (1982)
  153--195.

  \bibitem{Ginibre} J. Ginibre, G. Velo, The global Cauchy problem for the nonlinear Klein-Gordon
  equation, \emph{Math. Z.}, \textbf{189} (1985) 487--505.


  \bibitem{haraux87} A. Haraux, \emph{Semi-linear hyperbolic problems in bounded domains},
  Mathematical Reports, vol.~3, Harwood Academic, New York, 1987.





  \bibitem{lcc-tmna}C. Liu, F. Meng, Global well-posedness and attractor for damped wave equation
  with super-cubic nonlinearity and lower regular forcing on $\R^3$, \emph{Topol. Methods Nonlinear Anal.}, \textbf{49}(2) (2017) 551--563.

  \bibitem{lcc-jde} C. Liu, F. Meng, C. Sun, Well-posedness and attractors for a super-cubic
  weakly damped wave equation with $H^{-1}$ source term, \emph{J. Differential Equations}, \textbf{263} (2017) 8718--8748.

  \bibitem{lcc-aml} C. Liu, F. Meng, C. Zhong, A note on the global attractor for the weakly
  damped wave equation, \emph{Appl. Math. Lett.}, \textbf{41} (2015) 12--16.

  \bibitem{lcc26} C. Liu, F. Meng, Well-posedness and attractors for the sub-quintic weakly
  damped wave equation with low-regularity forcing on $\R^3$, to appear.

  \bibitem{masmoudi} N. Masmoudi, F. Planchon, On Uniqueness for the Critical Wave Equation,
  \emph{Communications in Partial Differential Equations}, \textbf{31} (2006) 1099--1107.

  \bibitem{mrw98} I. Moise, R. Rosa, X. Wang, Attractors for non-compact semigroups via
  energy equations, \emph{Nonlinearity}, \textbf{11} (1998) 1369--1393.

  \bibitem{planchon03} F. Planchon, On uniqueness for semilinear wave equations,
  \emph{Math. Z.}, \textbf{244} (2003) 587--599.

  \bibitem{ss} J. Shatah, M. Struwe, Well-Posedness in the Energy Space for Semilinear Wave
  Equations with Critical Growth, \emph{International Mathematics Research Notices},
  \textbf{7} (1994) 303--309.

  \bibitem{simon87} J. Simon, Compact sets in the space $L^p(0,T;B)$, \emph{Ann. Mat. Pura
    Appl.}, \textbf{146} (1987) 65--96.

  \bibitem{struwe99} M. Struwe, Uniqueness for critical nonlinear wave equations and wave
  maps via the energy inequality, \emph{Communications on Pure and Applied Mathematics},
  \textbf{52} (1999) 1179--1188.

  \bibitem{temam97} R. Temam, \emph{Infinite-dimensional dynamical systems in mechanics and
    physics}, 2nd ed., Applied Mathematical Sciences, vol.~68, Springer, New York, 1997.

  \bibitem{zelik15} S. Zelik, Strong uniform attractors for non-autonomous dissipative PDEs
  with non translation-compact external forces, \emph{Discrete Contin. Dyn. Syst. Ser. B},
  \textbf{20} (2015) 781--810.

\end{thebibliography}
\end{document}